\documentclass[a4paper]{gtart}

\usepackage[inner=20mm, outer=20mm, textheight=245mm]{geometry}

\usepackage[toc,page]{appendix}

\usepackage{psfrag, graphicx, subfigure, epsfig,color}

\usepackage{amsmath, amssymb, latexsym, euscript}

\usepackage[colorlinks,citecolor=blue,linkcolor=blue, backref=page]{hyperref}
\usepackage[nameinlink]{cleveref}

\usepackage[margin=1cm]{caption}

\usepackage{mathptmx}
\usepackage{wrapfig}
\usepackage[all,cmtip]{xy}
\usepackage{mathtools}
\usepackage{faktor}
\usepackage{tikz}
\usetikzlibrary{cd}
\usepackage{multicol}
\usepackage{listings}
\usepackage{stmaryrd} 

\usepackage[colorinlistoftodos]{todonotes}
\newcommand{\qt}{{/\!\!/}}

\newcommand{\Mod}[1]{\ (\mathrm{mod}\ #1)}

\newcommand{\Z}{\mathbb{Z}}
\newcommand{\C}{\mathbb{C}}

\newcommand{\N}{\mathbb{N}}
\newcommand{\F}{\mathbb{F}}
\newcommand{\K}{\mathbb{K}}

\DeclareMathOperator{\im}{im}

\DeclareMathOperator{\SL}{SL}

\DeclareMathOperator{\tr}{tr}

\DeclareMathOperator{\Hom}{Hom}

\DeclareMathOperator{\SLC}{\SL_2(\C)}
\DeclareMathOperator{\SLF}{\SL_2(\F)}

\def\X{X}
\def\R{R}

\def\P{\mathcal{P}}

\def\idMat{E}

\def\tree{\mathbf{T}}
\def\trFunc{I}

\def\glue{\Phi}
\def\klein{K}
\def\kleinGp{\Gamma_\klein}
\def\twistKlein{K\tilde{\times}I}

\def\seifert{\mathcal{S}}
\def\annulusM{\mathcal{R}}
\def\annulusa{\mathcal{A}}
\def\annulusb{\mathcal{B}}
\def\splitTorus{\mathcal{T}}

\def\glueSA{S_2}
\def\glueSB{S_3}
\def\glueRA{S_4}
\def\glueRApq{S_4(q,r)}
\def\glueRAp{\glueRA(q,1)}

\def\glueRB{S_5}

\def\curveI{(C1)}
\def\curveII{(C2)}
\def\curveIII{(C3)}
\def\surfI{(S1)}
\def\surfII{(S2)}
\def\surfIII{(S3)}
\def\surfIV{(S4)}
\def\surfV{(S5)}

\def\cTorusI{C_\glue}
\def\cTorusII{C_{\splitTorus,1}}
\def\cTorusIII{C_{\splitTorus,2}}

\def\cSA{C_{\glueSA}}
\def\cSB{C_{\glueSB}}

\def\cRA{C_\glue}
\def\cRB{C_{\glueRB}}
\def\csub{C}

\newcommand{\inv}[1]{#1^{-1}}
\newcommand{\bound}{\partial}
\newcommand{\pie}{\pi_1}

\theoremstyle{plain}
\newtheorem{theorem}{Theorem}
\newtheorem*{theorem*}{Theorem}
\newtheorem{lemma}[theorem]{Lemma}
\newtheorem{proposition}[theorem]{Proposition}
\newtheorem{corollary}[theorem]{Corollary}

\theoremstyle{definition}
\newtheorem{definition}[theorem]{Definition}
\newtheorem*{definition*}{Definition}

\newtheorem{conjecture}[theorem]{Conjecture}

\theoremstyle{remark}
\newtheorem{remark}[theorem]{Remark}

\numberwithin{equation}{section}

\usepackage{enumitem}
\newlist{enumRoman}{enumerate}{1}
\setlist[enumRoman]{label=(\Roman*)}

\begin{document}

\title{Essential tori in 3--manifolds not detected in any characteristic}
\author{Grace S. Garden, Benjamin Martin and Stephan Tillmann}

\begin{abstract} 
Infinite families of 3--dimensional closed graph manifolds and closed Seifert fibered spaces are exhibited, each member of which contains an essential torus not detected by ideal points of the variety of $\SL_2(\F)$--characters over any algebraically closed field $\F$. 
\end{abstract}


\makeshorttitle


\section{Introduction}
\label{sec:intro}

The ground-breaking work of Culler and Shalen \cite{Culler-Shalen-varieties-1983} detects essential surfaces in a 3--manifold by studying ideal points of curves in their $\SL_2(\C)$--character variety. Not all essential surfaces can be detected in this manner~(see \cite{Chesebro-boundary-2007,Schanuel-detection-2001}). Much work has been done to further characterise which surfaces are detected by which curves (for example, see \cite{Chesebro-closed-2013,Dunfield-incompressibility-2012,Segerman-detection-2011,Tillus-character-2004,Tillus-degenerations-2012, Yoshida-ideal-1991}). 
The theory also applies to so-called algebraic non-integral (ANI) representations, whose interplay with essential surfaces was crucial in the final resolution of the Smith conjecture~\cite{Morgan-Smith-1984}. These representations have also been shown not to detect all essential surfaces~\cite{Casella-ideal-2020}. 
More recently, the theory was generalised to
$\SL_n(\C)$--character varieties \cite{Hara-character-2021}, where it was shown that for every essential surface in a 3--manifold $M$, there is some $n$ such that it is detected by the $\SL_n(\C)$--character variety \cite{Friedl-representation-2018}. The theory is extended in another way in \cite{Garden-invitation-2024} by replacing the underlying field $\C$ with an arbitrary algebraically closed field $\F$ with characteristic $p$ and studying instead the variety of $\SL_2(\F)$--characters, which is closely related to the character variety. 

The motivating question for this paper is: given an essential surface $S$ in a compact 3--manifold $M$, is $S$ detected by an ideal point of the variety of $\SL_2(\F)$--characters of $M$ for some algebraically closed field $\F$? The answer to this question is unfortunately no for closed hyperbolic 3--manifolds. A closed, Haken, hyperbolic 3--manifold with no detected essential surface is given in \cite{Garden-invitation-2024}. However, essential spheres and tori have abelian fundamental group and may therefore be easier to detect because representations restricted to the surface are centralised by nontrivial 1--parameter groups. 

This paper analyses an infinite family of manifolds $N_\glue$ that consist of a twisted $I$--bundle over the Klein bottle (\Cref{sec:kleinbottle}) glued to the complement of the right-handed trefoil (\Cref{sec:trefoil}). These manifolds are either graph manifolds or Seifert fibered with base orbifolds $S^2(2,2,2,3)$ or $\mathbb{RP}^2(2,3)$. The essential surfaces in $N_\glue$ are classified (up to isotopy) in \Cref{sec:first_family}, and we determine which surfaces are detected in which characteristic in \Cref{sec:first_family}. Observations that follow from this analysis are:

\begin{theorem}\label{thm:graph-torus-never}
There are infinitely many graph manifolds with the property that they each contain (up to isotopy) a single connected essential surface, the single connected essential surface is a torus, and it is not detected by the variety of $\SL_2(\F)$--characters for any algebraically closed field $\F.$
\end{theorem}

\begin{theorem}\label{thm:graph-torus-char2}
There are infinitely many graph manifolds with the property that they each contain (up to isotopy) a single connected essential surface, the single connected essential surface is a torus, and it is detected by the variety of $\SL_2(\F)$--characters for an algebraically closed field $\F$ if and only if the characteristic of $\F$ is 2.
\end{theorem}

\begin{theorem}\label{thm:graph-torus-gen2}
There are infinitely many graph manifolds that contain (up to isotopy) two essential surfaces, a torus and a genus two surface, each of which is detected by the variety of $\SL_2(\F)$--characters for an algebraically closed field $\F$ if and only if the characteristic of $\F$ is 2.
\end{theorem}

\begin{theorem}\label{thm:graph-torus-nonsepgen2}
There are infinitely many graph manifolds that contain (up to isotopy) two essential surfaces, a torus and a non-separating genus two surface, and with:
\begin{itemize}
	\item the torus not detected by the variety of $\SL_2(\F)$--characters for any algebraically closed field $\F$, and 
	\item the non-separating genus two surface detected by the variety of $\SL_2(\F)$--characters for every algebraically closed field $\F$.
\end{itemize}
\end{theorem}

The second part of the above theorem follows because each non-separating orientable surface is Poincar\'e dual to an epimorphism from the fundamental group of the ambient 3--manifold to the integers, and hence detected by a curve of reducible representations. 

\begin{theorem}\label{thm:RP2family}
There are infinitely many Seifert fibered 3--manifolds with base $\mathbb{RP}^2(2,3)$ and that contain (up to isotopy) exactly two connected essential surfaces, both vertical tori, with:
\begin{itemize}
	\item one torus is not detected by the variety of $\SL_2(\F)$--characters for any algebraically closed field $\F$, and 
	\item the other torus is detected by the variety of $\SL_2(\F)$--characters for an algebraically closed field $\F$ if and only if the characteristic of $\F$ is 2.
\end{itemize}
\end{theorem}

\begin{theorem}\label{thm:S2family}
There are infinitely many Seifert fibered 3--manifolds with base $S^2(2,2,2,3)$ and that contain infinitely many pairwise non-isotopic essential tori, all of which are detected by the variety of $\SL_2(\F)$--characters for every algebraically closed field $\F$. 
\end{theorem}
For \Cref{thm:S2family}, we remark that the theorem is not proved for \emph{all} essential tori in the stated Seifert fibered 3--manifolds, but that we only prove it for an infinite family thereof. 

\Cref{thm:graph-torus-char2,thm:graph-torus-gen2,thm:RP2family} give examples of manifolds with an essential surface that is not detected by the $\SLC$--character variety but is detected by the variety of $\SLF$--characters for $p=2$; 
\Cref{thm:graph-torus-never,thm:graph-torus-nonsepgen2,thm:RP2family} give examples of manifolds with an essential surface 
 that is not detected by the variety of $\SLF$--characters for all $p$.

\textbf{Acknowledgements.} 
Research of the first author is supported by an Australian Government Research Training Program scholarship.
The second author thanks Xingru Zhang for useful discussions and the School of Mathematics and Physics at the University of Queensland for hospitality.
Research of the third author is supported in part under the Australian Research Council's ARC Future Fellowship FT170100316. The authors thank Eric Chesebro and Daniel Mathews for helpful comments on an earlier draft of the manuscript.

For the purpose of open access, the authors have applied a Creative Commons Attribution (CC BY) licence to any Author Accepted Manuscript version arising from this submission.


\section{Preliminaries}
\label{sec:preliminaries}


For the remainder of the paper, let $\F$ be an algebraically closed field of characteristic $p\geq 0.$ 
We review the extension of Culler--Shalen theory presented in \cite{Garden-invitation-2024}. This outlines how to detect essential surfaces in a 3--manifold by studying ideal points of curves in the variety of $\SL_2(\F)$--characters.

We write $\F = \F_p$ if we want to emphasise the characteristic and $\Z_p$ for the finite field of $p$ elements. 

\subsection{The variety of characters}

Let $M$ be an orientable, compact $3$--manifold. The variety of $\SLF$--characters of $M$ is written $X(M,\F)$. 
We use $E$ to denote the identity matrix in $\SLF$. See \cite{Garden-invitation-2024} for a detailed discussion of $X(M,\F);$ this paper only requires a few facts that are familiar from the standard material~\cite{Shalen-representations-2002} which we now summarise.

The \textbf{$\SLF$--representation variety} of $\Gamma$ is the space of representations $\rho\co \Gamma \to \SLF$,
\begin{equation}
	\R(\Gamma,\F) = \Hom(\Gamma,\SLF).
\end{equation}

Given a finite, ordered generating set $( \gamma_1, \ldots, \gamma_n)$ of $\Gamma,$ the Hilbert basis theorem implies that we can imbue $\R(\Gamma,\F)$ with the structure of an affine algebraic subset of $\SLF^n \subset \F^{4n}$ via the natural embedding
\begin{align*}
	\varphi\co \R(\Gamma,\F) &\hookrightarrow \SLF^n \subset \F^{4n}\\
	\rho &\mapsto (\rho(\gamma_1), \ldots, \rho(\gamma_n)).
\end{align*}
In particular, we identify $R(F_n,\F)$ with $\SLF^n$, where $F_n$ is the free group on generators $\gamma_1,\ldots, \gamma_n$.
The Zariski topology on $\F^{4n}$ induces a topology on $\R(\Gamma,\F).$

The group $\SLF$ acts by conjugation on $\SLF^n.$ The action preserves the natural embedding of $\R(\Gamma,\F)$ and is algebraic.  Since $\SLF$ is reductive, we can form the quotient variety $R(\Gamma,\F)\qt \SLF$.  The co-ordinate ring $\F[R(\Gamma,\F)\qt \SLF]$ is the ring of invariants $\F[R(\Gamma,\F)]^{\SLF}$.

Given $\gamma\in \Gamma$, define $\trFunc_\gamma \in \F[R(\Gamma,\F)]$ by 
\[
\begin{split}
	\trFunc_\gamma\co R(\Gamma,\F) 	&\to \F\\
			\rho			&\mapsto \tr(\rho(\gamma)).
\end{split}
\]
Let ${\mathcal T}_\Gamma$ be the subring of $\F[R(\Gamma,\F)]^{\SLF}$ generated by the $I_\gamma$.  It can be shown that $\F[R(\Gamma,\F)]^{\SLF}$ is finite as a ${\mathcal T}_\Gamma$-module; in particular, ${\mathcal T}_\Gamma$ is finitely generated as a $k$-algebra \cite[Theorem 1.4]{Martin-2003}.  Hence there is an affine variety $X(\Gamma,\F)$ whose co-ordinate ring is ${\mathcal T}_\Gamma$.  The inclusion of ${\mathcal T}_\Gamma$ in $\F[R(\Gamma,\F)]^{\SLF}$ gives rise to a morphism of varieties $q_\Gamma\colon R(\Gamma,\F)$ to $X(\Gamma,\F)$.  If $p= 0$ or if $\Gamma= F_n$ for some $n\in \N$ then ${\mathcal T}_\Gamma= \F[R(\Gamma,\F)]^{\SLF}$, so in this case we may identify $X(\Gamma,\F)$ with $R(\Gamma,\F)\qt \SLF$.  For more details and further discussion, see \cite{Garden-invitation-2024} and \cite{Martin-2024}.

We say that $\rho, \sigma \in \R(\Gamma,\F)$ are \textbf{closure-equivalent} if the Zariski closures of their orbits intersect; we write $\rho\sim_c \sigma$.
The \textbf{character of a representation $\rho$} is the trace function
\begin{equation}
\begin{split}
	\tau_\rho\co \Gamma &\to \F \\
	\gamma&\mapsto\tr(\rho(\gamma)).
\end{split}	
\end{equation}

The next result is \cite[Corollary 16]{Garden-invitation-2024}.
\begin{corollary}
\label{cor:trace_function_characterisation}
Let $\Gamma$ be a finitely generated group. 
Suppose $\rho, \sigma \in \R(\Gamma,\F).$ Then the following four statements are equivalent:
\begin{enumerate}
	\item $\tau_\rho = \tau_\sigma$;
	\item $\tau_\rho$ and $\tau_\sigma$ agree on all ordered products of distinct generators;
	\item $\tau_\rho$ and $\tau_\sigma$ agree on all ordered single, double and triple products of distinct generators;
	\item $\rho$ and $\sigma$ are closure-equivalent.
\end{enumerate}
\end{corollary}

\noindent It follows from Corollary~\ref{cor:trace_function_characterisation} that we may identify $\X(\Gamma,\F)$ with the quotient space of the equivalence relation $\sim_c$.  The map $q_\Gamma$ factorises as $j_\Gamma\circ \pi_\Gamma$, where $\pi_\Gamma\colon \R(\Gamma,\F)\to \R(\Gamma,\F)\qt \SLF$ is the canonical projection and $j_\Gamma\colon \R(\Gamma,\F)\qt \SLF\to \X(\Gamma,\SLF)$ is a morphism.  It follows from the above discussion that $j_\Gamma$ is finite and bijective, so $j_\Gamma$ is a homeomorphism.  By standard geometric invariant theory, $\pi_\Gamma$ maps closed $\SLF$--stable subsets of $\R(\Gamma,\F)$ to closed subsets of $\R(\Gamma,\F)\qt \SLF$, so $j_\Gamma$ maps closed $\SLF$--stable subsets of $\R(\Gamma,\F)$ to closed subsets of $\X(\Gamma,\F)$.


We say that a subgroup of $\SLF$ is \textbf{reducible} if its action on $\F^2$ has an invariant 1--dimensional subspace. Otherwise it is \textbf{irreducible}. 
We call a \emph{representation} (or equivalently an $n$--tuple of matrices) \textbf{irreducible} (resp. \textbf{reducible}) if it generates an {irreducible} (resp. {reducible}) subgroup of $\SLF$. 
We say that a \emph{character} is \textbf{irreducible} (resp. \textbf{reducible}) if it is associated with an irreducible (resp. {reducible}) $n$--tuple of matrices. 

In particular, all abelian subgroups of $\SLF$ are reducible. We call a representation (or equivalently an $n$--tuple of matrices) \textbf{abelian} if it generates an {abelian} subgroup of $\SLF.$

We define
$$ \R^{\rm red}(\Gamma,\F)= \{\rho\in \R(\Gamma,\F)\,|\,\mbox{$\rho$ is reducible}\}.  $$
Define $\R^{\rm irr}(\Gamma,\F)$ to be the closure of the complement of $\R^{\rm red}(\Gamma,\F)$ in $\R(\Gamma,\SLF)$. Similarly, define $\X^{\rm red}(\Gamma,\F)= q_\Gamma(\R^{\rm red}(\Gamma,\F))$ and $\X^{\rm irr}(\Gamma,\F)= q_\Gamma(\R^{\rm irr}(\Gamma,\F))$.  

The results above imply that both $\X^{\rm red}(\Gamma,\F)$ and $\X^{\rm irr}(\Gamma,\F)$ are closed.  One can show using the discussion in \cite[\S 8]{Martin-2003} that if $\rho\in \R(\Gamma,\F)$ is irreducible then $j_\Gamma(\rho)\not\in \X^{\rm red}(\Gamma,\F)$; hence $\X^{\rm irr}(\Gamma,\F)$ is the closure of the complement of $\X^{\rm red}(\Gamma,\F)$ in $\X(\Gamma,\SLF)$.

We will use the following result \cite[Lemma 7]{Garden-invitation-2024} repeatedly.
\begin{lemma}
\label{lem:irred_crit}
 Let $A,B\in \SLF$.  Then $A$ and $B$ generate a reducible subgroup of $\SLF$ if and only if $\tr[A,B]= 2$.
\end{lemma}

If $\Gamma = \pi_1(M)$ for a compact manifold $M$, then we also write $\R(M,\F) = \R(\Gamma,\F)$ and $\X(M,\F) = \X(\Gamma,\F)$; likewise for the respective irreducible and reducible components. 


\subsection{Essential surfaces}

We define essential surfaces following standard terminology from Jaco \cite{Jaco-lectures-1980}.
A {surface} $S$ in a compact 3--manifold $M$ will always mean a 2--dimensional piecewise linear submanifold \emph{properly embedded} in $M$. That is, a closed subset of $M$ with $\partial S = S \cap \partial M$. If $M$ is not compact, we replace it by a compact core. 

An embedded sphere $S^2$ in a 3--manifold $M$ is called \textbf{incompressible} or \textbf{essential} if it does not bound an embedded ball in $M$, and a 3--manifold is \textbf{irreducible} if it contains no incompressible 2--spheres. 
An orientable surface $S$ without 2--sphere or disc components in an orientable 3--manifold $M$ is called \textbf{incompressible} if for each disc $D\subset M$ with $D \cap S = \partial D$ there is a disc $D' \subset S$ with $\partial D' = \partial D$.
Note that \emph{every} properly embedded disc in $M$ is incompressible.

A surface $S$ in a 3--manifold $M$ is \textbf{$\partial$--compressible} if either 
\begin{enumerate}
\item $S$ is a disc and $S$ is parallel to a disc in $\partial M,$ or
\item $S$ is not a disc and there exists a disc $D\subset M$ such that $D\cap S = \alpha$ is an arc in $\partial D,$ $D\cap \partial M = \beta$ is and arc in $\partial D,$ with $\alpha \cap \beta = \partial \alpha = \partial \beta$ and $\alpha \cup \beta = \partial D$ and either $\alpha$ does not  separate $S$ or $\alpha$ separates $S$ into two components and the closure of neither is a disc.
\end{enumerate}
Otherwise $S$ is \textbf{$\partial$--incompressible}.

\begin{definition}  \cite{Shalen-representations-2002}\label{def:essential}
A surface $S$ in a compact, irreducible, orientable 3--manifold is said to be \textbf{essential} if it has the following properties:
    \begin{enumerate}
       \item $S$ is bicollared;
       \item the inclusion homomorphism $\pi_1(S_i) \to \pi_1(M)$ is
          injective for
           every component $S_i$ of $S$;
       \item no component of $S$ is a 2--sphere;
       \item no component of $S$ is boundary parallel;
       \item $S$ is non-empty.
    \end{enumerate}
\end{definition}
In the first condition, bicollared means $S$ admits a map $h\co S\times [-1,1] \to M$ that is a homeomorphism onto a neighbourhood of $S$ in $M$ such that $h(x,0)=x$ for every $x \in S$ and $h(S\times [-1,1])\cap \bound M = h(\bound S\times [-1,1])$.
The surface $S$ being bicollared in orientable $M$ implies $S$ is orientable. 
The second condition is equivalent to saying that there are no compression discs for the surface (cf.~\cite[Lemma 6.1]{Hempel-manifolds-1976}). Hatcher~\cite[Lemma 1.10]{Hatcher-notes} implies that if each boundary component of the essential surface $S$ lies on a torus boundary component of $M$, then $S$ is both incompressible and $\partial$--incompressible.

A compact, irreducible 3--manifold that contains an essential surface is called \textbf{Haken}.


Essential surfaces in Seifert fibered manifolds can be described more explicitly. In any connected, compact, irreducible Seifert fibered manifold $M$, Hatcher shows \cite[Proposition 1.11]{Hatcher-notes} that each essential surface is isotopic to either a \textbf{vertical} surface (a union of regular fibres) or a \textbf{horizontal} surface (transverse to the fibration, giving a branched covering of the base orbifold with branch points corresponding to the intersections with singular fibres).

Let the base orbifold of $M$ be $Q$ and let $\hat{Q}$ denote the surface obtained by removing regular neighbourhoods of cone points in $Q$. By a result of Schultens, ~\cite[Lemma 39]{Schultens-kakimizu-2018}, vertical essential surfaces are in bijective correspondence with the isotopy classes of simple closed curves in $\hat{Q}$.


\subsection{Essential surfaces detected by ideal points}


Let $C$ be a closed irreducible curve in $\F^m$ for some $m$.  We define the \textbf{projective completion} $C'$ of $C$ to be the closure of the image of $C$ under the map $J: \F^m\to \F P^m$ defined by $J(z_1,\ldots,z_m)=[1,z_1,\ldots,z_m]$.  We define $\tilde{C}$ to be the normalisation of the projective completion of $C$; then $\tilde{C}$ is a smooth irreducible projective curve.  Let $f\colon \tilde{C}\to C'$ be the normalisation map.  The \textbf{ideal points} of $C$ are the points $\xi\in \tilde{C}$ such that $f(\xi)\in \overline{J(C)}-J(C)$.  It can be shown that $\tilde{C}$ and the notion of an ideal point do not depend on the choice of embedding of $C$ in $\F^m$.  Note that the function fields $\F(\tilde{C})$ and $\F(C)$ are isomorphic.  

Below we consider closed irreducible curves $C$ in irreducible components of $\X(M, \F).$  Of course, if an irreducible component of $\X(M, \F)$ has dimension 0 (resp., dimension 1) then there are no such curves (resp., exactly one such curve).

Take an ideal point $\xi\in\tilde{C}$. This defines a \textbf{discrete rank 1 valuation} on $\F(C)$ via
\begin{equation}
	\begin{split}
	v_\xi\co \F(C)&\to\Z\cup\{\infty\}\\
	v_\xi(f)&=
		\begin{cases}
		q & f\text{ has a zero of order }q\text{ at }\xi,\\
		\infty & f=0, \\
		-q & f\text{ has a pole of order }q\text{ at }\xi.
		\end{cases}
	\end{split}
\end{equation}
Write $\Gamma=\pi_1(M)$.
Following Culler and Shalen, we assume that $C$ maps to a curve in $\X(M, \F).$ Further, we assume that $\xi$ maps to an ideal point of that curve. 
This implies that there is an element $\gamma \in \Gamma$ such that $v_\xi(\trFunc_\gamma) < 0.$ We can now use the valuation to construct the {Bass-Serre tree} $\tree_\xi$ with an action of $\SL_2(\F(C))$. We can then use the \textbf{tautological representation}
\begin{equation}
\begin{split}
	\P\co \pie(M)&\to \SL_2(\F(C)),\\
	\P(\gamma)&=\begin{pmatrix} a & b\\ c& d \end{pmatrix} \text{ for }\rho(\gamma)=\begin{pmatrix} a(\rho) & b(\rho)\\ c(\rho)& d(\rho) \end{pmatrix}
\end{split}
\end{equation}
to get an action of $\pi_1(M)$ on $\tree_\xi$. The following is \cite[Property 5.4.2]{Shalen-representations-2002}, note the proof applies verbatim.
\begin{proposition}
\label{pro:stabilisers}
Let $C$ be a curve in $\X(\Gamma, \F).$ To each ideal point $\xi$ of $\tilde{C},$ one can associate a splitting of $\Gamma$ with the property that for each element $\gamma \in \Gamma,$ the following are equivalent:
\begin{enumerate}
\item $v_\xi(\trFunc_\gamma)\ge 0$
\item A vertex of the Bass-Serre tree $\tree_\xi$ is fixed by $\gamma.$
\end{enumerate}
\end{proposition}

A construction by Stallings~\cite{Stallings-topological-1965} can now be used to associate an essential surface $S$ to the ideal point $\xi$ (see \cite{Shalen-representations-2002} for details). The surface is not unique and we usually identify parallel components. If a given essential surface $S$ can be associated to an ideal point $\xi,$ then we say that \textbf{$S$ is detected by $\xi.$} 
The components of a detected surface $S$ and the complementary components $M \setminus S$ satisfy a number of key properties. For instance, \Cref{pro:stabilisers} directly implies:

\begin{corollary}
\label{cor:stabilisers_cor}
If $S$ is detected by the ideal point $\xi$, then
\begin{enumerate}[label=(\roman*)]
	\item for each component $M_i$ of $M-S$, the subgroup $\text{im}(\pi_1(M_i)\to \pi_1(M))$ of $\pi_1(M)$ is contained in the stabiliser of a vertex of $\tree_\xi$; and 
	\item for each component $S_j$ of $S$, the subgroup $\text{im}(\pi_1(S_j)\to \pi_1(M))$ of $\pi_1(M)$ is contained in the stabiliser of an edge of $\tree_\xi$.
\end{enumerate}
\end{corollary}

A \textbf{boundary slope} is the slope of the boundary curve of an essential surface $S$ in $M$ that has non-empty intersection with $\bound M$. 
The relationship between boundary slopes of essential surfaces in $M$ and valuations of trace functions of peripheral elements is given in the following reformulation of \cite[Proposition 1.3.9]{Culler-Dehn-1987}.

\begin{proposition}
\label{pro:detected_surface}
Let $M$ be a compact, orientable, irreducible 3--manifold with $\partial M$ a torus. 
Let $C$ be a curve in the variety of characters with an ideal point $\xi$.
We have the following mutually exclusive cases.
\begin{enumerate}
  \item If there is an element $\gamma$ in $\im (\pi_1(\partial M) \to \pi_1(M))$
     such that $v_\xi(\trFunc_\gamma)<0$, then up to inversion there
     is a unique primitive element $\alpha \in \im(\pi_1(\partial M)\to \pi_1(M))$ such that 
     $v_\xi(\trFunc_\alpha)\ge 0$. Then every essential surface $S$ detected by $\xi$ has non-empty boundary and 
     $\alpha$ is parallel to its boundary
     components.  
  \item If $v_\xi(\trFunc_\gamma)\ge 0$ for all 
     $\gamma \in \im (\pi_1(\partial M)\to\pi_1(M))$, 
     then an essential surface $S$ detected by $\xi$ may be chosen that is disjoint from $\partial M$.
\end{enumerate}
\end{proposition}


\section{The Klein bottle group}
\label{sec:kleinbottle}

We consider the Klein bottle group via a well-known presentation,
\begin{equation}
\label{eqn:kleinbottle}
	\kleinGp= \langle a, b \mid aba^{-1}b = 1 \rangle.
\end{equation}

\begin{figure}[bth]
    \begin{center}
     \includegraphics[height=4.4cm]{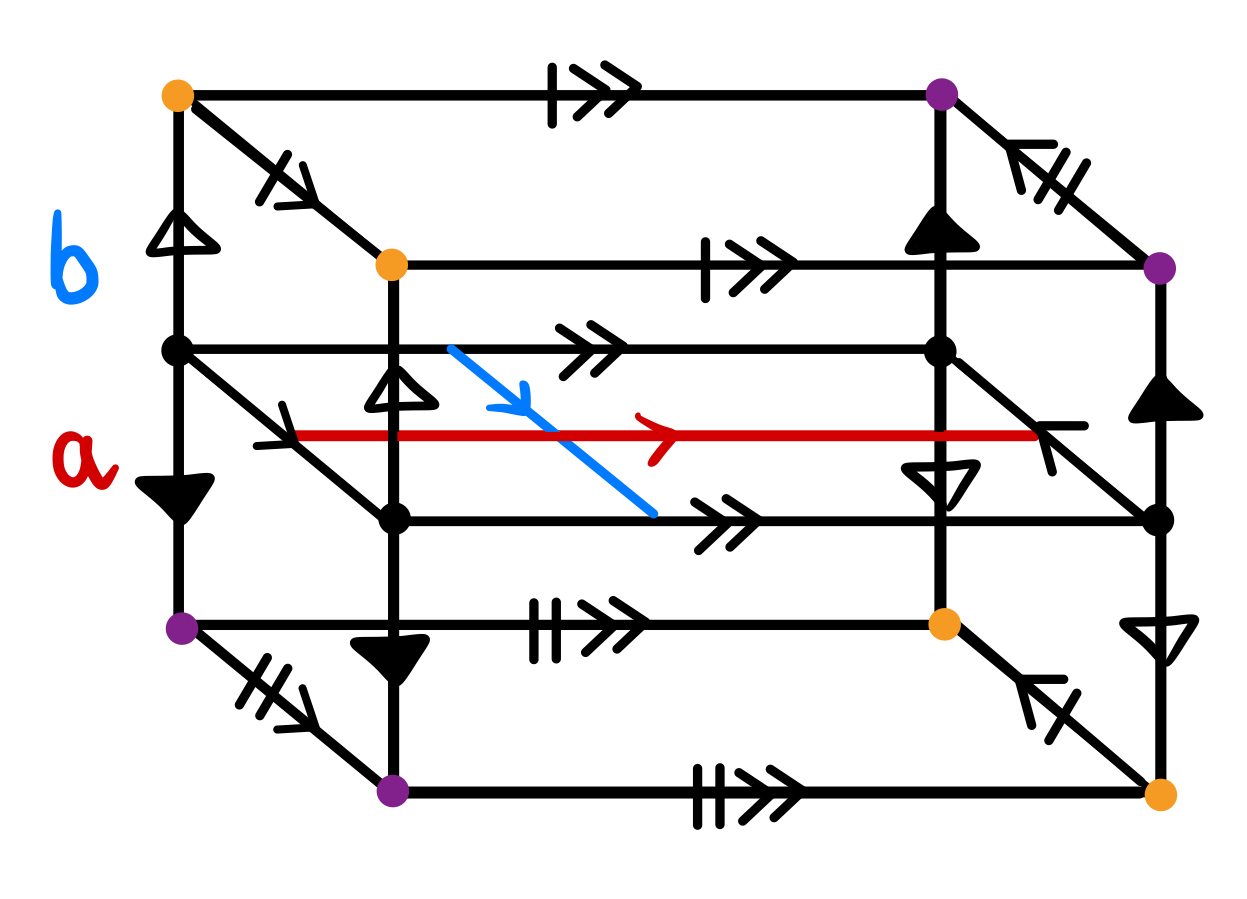}
          \end{center}
    \caption{Generators for the twisted $I$--bundle over the Klein bottle.}
    \label{fig:Klein1}
\end{figure}

The following result is proven in \cite{Friedl-linear-2018}.

\begin{proposition}[\cite{Friedl-linear-2018}]
\label{prop:kleinbottle}
Let $\K$ be a field and $\kleinGp$ the Klein bottle group as in \Cref{eqn:kleinbottle}. Suppose 
 $\sigma \colon \kleinGp \to \SL_2(\K)$ is a representation. If $\sigma(b) \neq \pm\idMat$, then one of the following occurs:
\begin{enumerate}
	\item[(i)] $\sigma$ factors through the abelianization of $\kleinGp$.
	\item[(ii)] $\sigma(a^2) = \pm\idMat$.
\end{enumerate}
In particular, $\kleinGp$ does not admit faithful representations into $\SL_2(\K)$ for any field $\K$.
\end{proposition}
Note, $\sigma(a^2)=\sigma(a)^2=\pm\idMat$ implies either $\sigma(a)=\pm\idMat$ or $\tr(\sigma(a))=0$. 

Consider the variety of characters $\X(\klein,\F)$. We can easily compute the subvarieties $\X^{\text{irr}}(\klein,\F)$ and $\X^{\text{red}}(\klein,\F)$ respectively. The main observation is that in characteristic 2 there are infinitely many elements in $\SLF$ with order two and, as a result, we find infinitely many non-conjugate abelian representations for $\kleinGp$.

Each irreducible representation is conjugate to the form
\begin{equation}
\label{eqn:irredRepnKlein}
	\rho(a) = \begin{pmatrix} 0 & -1 \\ 1 & 0\end{pmatrix}, \quad \rho(b) = \begin{pmatrix} y & 0 \\ 0 & y^{-1}\end{pmatrix}, \quad y\in\F\setminus\{0\}.
\end{equation}
By~\Cref{lem:irred_crit}, this representation is indeed irreducible if and only if 
\[
0 \neq \tr [\rho(a) , \rho(b) ]  - 2 = y^2 - 2 + y^{-2}.
\]
Hence the representation is irreducible if and only if $y^2 \neq 1.$
Letting $t = y+y^{-1},$ the map $$\rho \mapsto \left(\tr(\rho(a)),\tr(\rho(b)),\tr(\rho(ab))\right) \in \F^3$$ gives
\begin{equation*}
	\X^{\text{irr}}(\klein,\F) = \{(0,t,0)\mid t\in\F\}.
\end{equation*}
From the trace condition on the commutator, we know $\X^{\text{irr}}(\klein,\F)\cap\X^{\text{red}}(\klein,\F)=\{(0,\pm2,0)\}$. We now determine all reducible representations up to conjugacy. A reducible representation is conjugate to either
\begin{align}
	\label{eqn:redRepnKlein1} 
	\rho(a) &= \begin{pmatrix} x & 1 \\ 0 & x^{-1}\end{pmatrix}, \quad \rho(b) = \begin{pmatrix} \pm1 & 0 \\ 0 & \pm1\end{pmatrix}, \quad x\in\F\setminus\{0\}, \text{ or}\\
	\label{eqn:redRepnKlein2} 
	\rho(a) &= \begin{pmatrix} x & u \\ 0 & -x\end{pmatrix}, \quad \rho(b) = \begin{pmatrix} \pm1 & 1 \\ 0 & \pm1\end{pmatrix}, \quad u\in\F,\: x\in\F\setminus\{0\},\: x^2=-1.
\end{align}
Here we have $\rho(b)=\pm\idMat$ for the representations in  \Cref{eqn:redRepnKlein1} and $\tr(\rho(a))=0$ for the representations in \Cref{eqn:redRepnKlein2}.

In characteristic $p=2$,
\[
	\X^{\text{red}}(\klein,\F) = \{ (s,0,s)\mid s\in\F\}
\]
and the matrix pair in \Cref{eqn:redRepnKlein2} splits into infinitely many conjugacy classes, 
\begin{align}
	\label{eqn:redRepnKleinChar21}
	\rho(a) = \begin{pmatrix} 1 & u \\ 0 & 1\end{pmatrix}, \quad \rho(b) = \begin{pmatrix} 1 & 1 \\ 0 & 1\end{pmatrix}, \quad u\in\F.
\end{align}

In characteristic $p\neq2$,
\[
	\X^{\text{red}}(\klein,\F) = \{ (s,2,s)\mid s\in\F\}\cup\{ (s,-2,-s)\mid s\in\F\}
\]
and the matrix pair in \Cref{eqn:redRepnKlein2} is conjugate to
\begin{equation}
\label{eqn:redRepnKleinCharNon22}
	\hspace{-0.56cm}\rho(a) = \begin{pmatrix} x & 0 \\ 0 & -x\end{pmatrix}, \quad \rho(b) = \begin{pmatrix} \pm1 & 1 \\ 0 & \pm1\end{pmatrix}, \quad x\in\F\setminus\{0\},\: x^2=-1.
\end{equation}

Let ${\twistKlein}$ be a twisted $I$--bundle over the Klein bottle with $\Gamma_{\twistKlein}$ the associated fundamental group,
\[
	\Gamma_{\twistKlein} = \kleinGp = \left \langle a,b  \mid  ab\inv{a}b \right\rangle.
\]

Then the variety of characters of $\twistKlein$ is the same as that of $\klein$, 
\[
\begin{split}
	\X^{\text{irr}}(\twistKlein,\F)	&=\X^{\text{irr}}(\klein,\F),\\
	\X^{\text{red}}(\twistKlein,\F)	&=\X^{\text{red}}(\klein,\F).
\end{split}
\]
We note from the above calculations that in any characteristic $p\ge 0$, the variety $\X^{\text{irr}}(\twistKlein,\F)$ is a curve with a single ideal point, and $\X^{\text{red}}(\twistKlein,\F)$ consists of either a single curve with a single ideal point (if $p=2$) or two disjoint curves, each with a single ideal point (if $p\neq 2$).

The peripheral subgroup for $\Gamma_{\twistKlein}$ is 
\[
	P_{\twistKlein}=\langle a^2, b\rangle.
\]

The twisted $I$--bundle over the Klein bottle has two Seifert fibrations. See \cite[\S1.5]{Aschenbrenner-book-2015} for a succinct discussion; we summarise the points needed for our application and use the above description of the variety of characters to determine which components of $\X(\twistKlein,\F)$ detect which essential surface.

The two Seifert fibrations of ${\twistKlein}$ have base orbifolds $D^2(2,2)$ and the M\"obius band (with no cone points) respectively. 
We know there are two types of connected surfaces in $\twistKlein$: vertical and horizontal. 
The fibration over the M\"obius band shows that there is one connected horizontal surface that arises from an unbranched double cover of the M\"obius band (a separating annulus $\annulusa$ with boundary slope $a^2$) and one connected vertical surface (a non-separating annulus $\annulusb$ with boundary slope $b$). 
Note that these surfaces cannot be isotoped to be disjoint. It follows that an essential surface in ${\twistKlein}$ consists either of parallel copies of $\annulusa$ or of parallel copies of $\annulusb$. Moreover, an essential surface in ${\twistKlein}$ is uniquely determined by its boundary curves.

The slope $a^2$ is a regular fibre in the Seifert fibration of ${\twistKlein}$ with base orbifold $D^2(2,2)$. The annulus $\annulusa$ is foliated by regular fibres and decomposes ${\twistKlein}$ into two solid tori with Seifert structures over the base orbifold $D^2(2).$ Hence, $\annulusa$ is an essential annulus in the boundary of the two solid tori. The presentation of the fundamental group corresponding to this is a free product with amalgamation of the form
\[
	\langle c \rangle \star_{c^2 = d^2} \langle d \rangle.
\]
Here, $c^2=d^2$ corresponds to the regular fibre. An isomorphism to the presentation given in \Cref{eqn:kleinbottle} is:
\[
	c\mapsto a \qquad\text{and}\qquad d\mapsto b^{-1}a.
\]
Note that this gives $c^2 = d^2 \mapsto a^2$ as claimed.

For any characteristic $p\ge 0$, the restriction of $\trFunc_a$ (and hence of  $\trFunc_{a^2}$) to the curve $$\X^{\text{irr}}(\twistKlein,\F)=\{(0,t,0)\mid t\in\F\}$$ is constant equal to zero while the trace function $\trFunc_b$ has a pole of order one at the ideal point. It follows from \Cref{pro:detected_surface} and the above classification of essential surfaces that $\annulusa$ is detected by the curve $\X^{\text{irr}}(\twistKlein,\F)$ and it also follows that $\annulusb$ is not detected by this curve. Also note that the fundamental groups of the complementary solid tori are generated by $\langle a \rangle$ and $\langle b^{-1}a \rangle$ respectively. Now the trace functions $\trFunc_{a}$ and $\trFunc_{b^{-1}a}$ are constant equal to zero on $\X^{\text{irr}}(\twistKlein,\F),$ thus verifying the splitting detected at the ideal point.

\begin{figure}[t]
    \begin{center}
    \subfigure[Image of $\annulusa$ in $D^2(2,2)$]{%
        \label{fig:KleinD22}%
        \includegraphics[height=4.5cm]{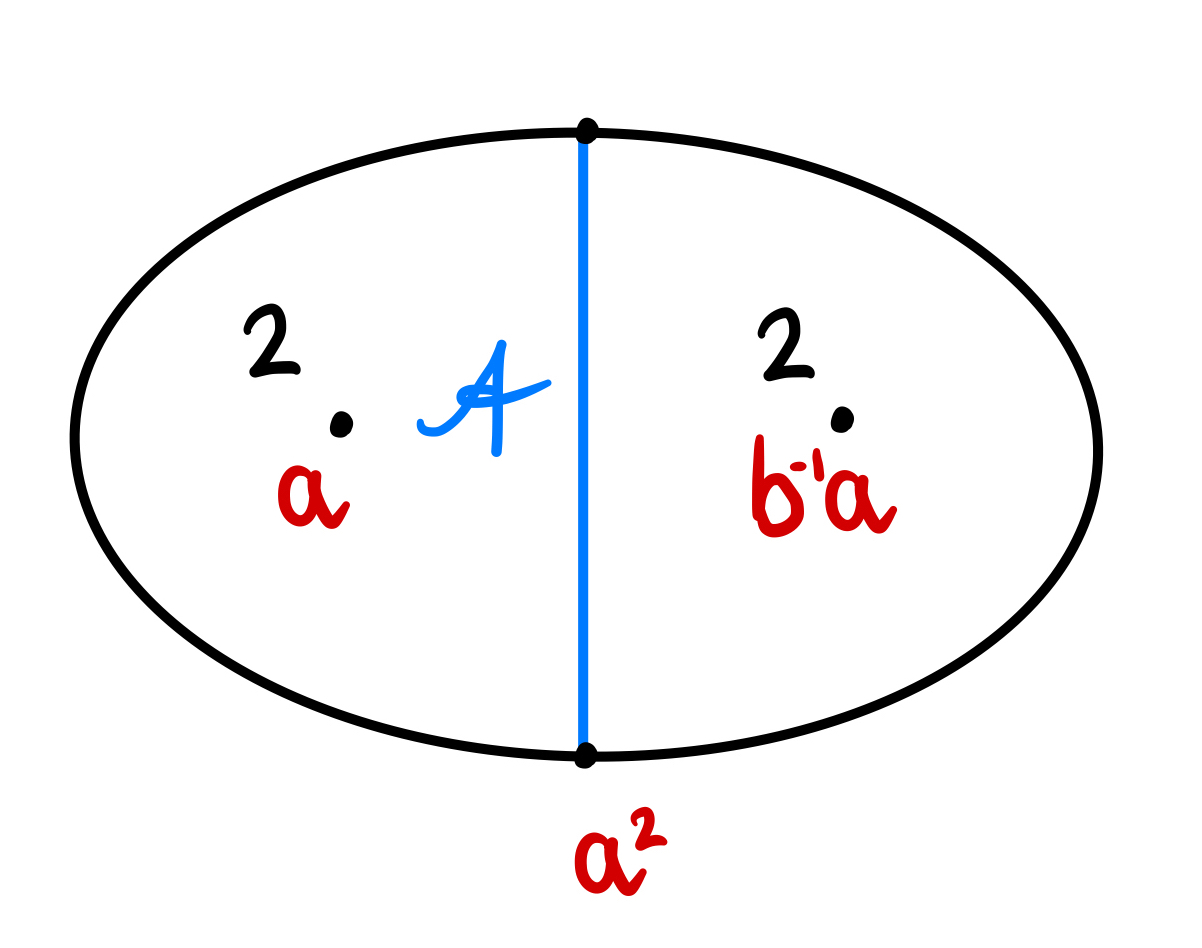}}
    \qquad
     \subfigure[Image of $\annulusb$ in the M\"obius band]{%
        \label{fig:KleinMobius}%
        \includegraphics[height=4.5cm]{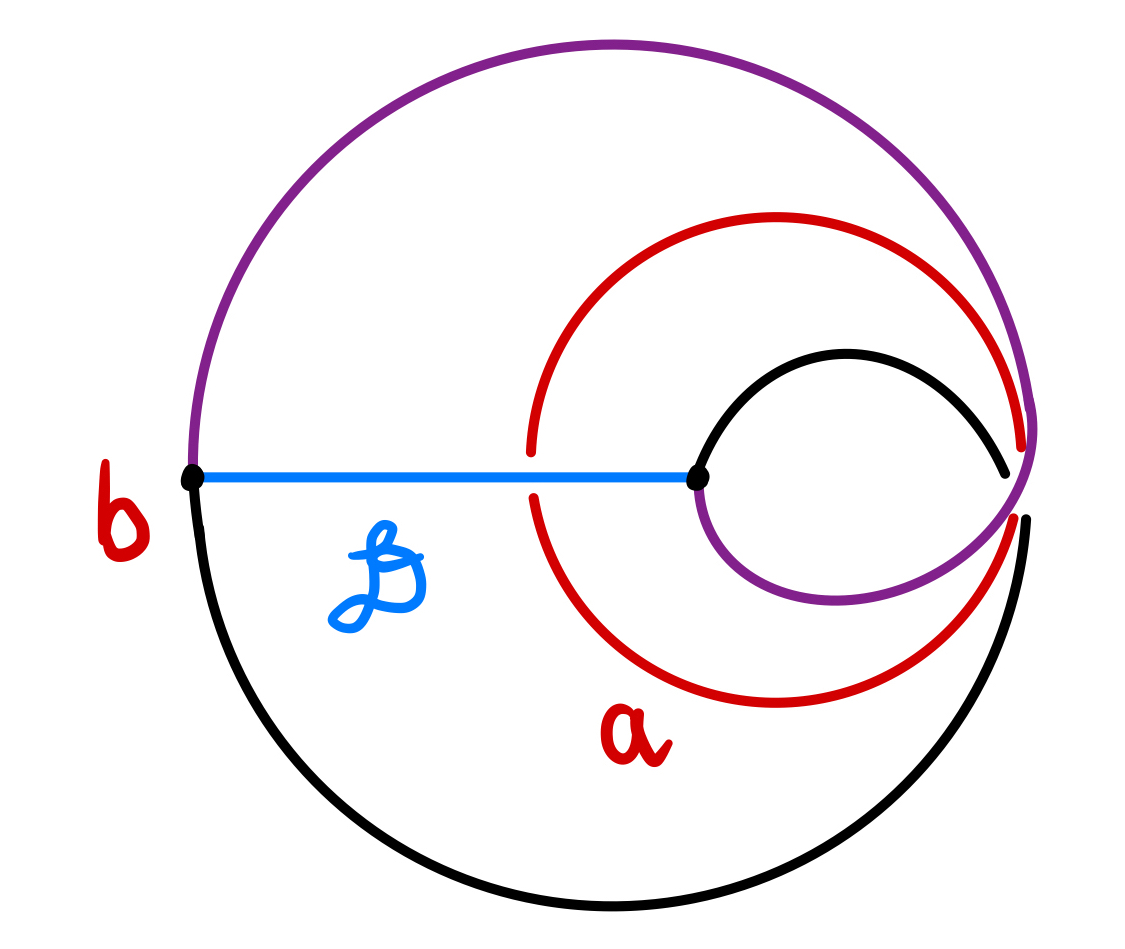}}
          \end{center}
    \caption{Essential annuli in the twisted $I$--bundle over the Klein bottle.}
    \label{fig:Klein}
\end{figure}

The slope $b$ is a regular fibre in a Seifert fibration with base orbifold the M\"obius band (with no cone points). 
The non-separating annulus $\annulusb$ is foliated by regular fibres in this fibration, hence $\im\left(\pi_1(\annulusb) \to \pi_1(\twistKlein)\right) = \langle b \rangle.$ Now the trace function $\trFunc_{b}$ is constant on each component of $\X^{\text{red}}(\twistKlein,\F)$ (there is one curve for $p=2$ and two curves for $p\neq 2$) and the trace function $\trFunc_{a}$ has a pole. In particular, \Cref{pro:detected_surface} implies that $\annulusb$ is detected by each curve in $\X^{\text{red}}(\twistKlein,\F)$ and  $\annulusa$ is not detected by the curves in $\X^{\text{red}}(\twistKlein,\F)$. Note that this is related to the fact that $\annulusb$ is dual to the epimorphism $\pi_1(\twistKlein) \to \Z$ defined by $a\mapsto 1$ and $b\mapsto 0.$ 


\section{The right-handed trefoil group}
\label{sec:trefoil}

We consider the complement of the right-handed trefoil $M$, and the associated fundamental group via a well-known presentation, 
\[
	\Gamma_M = \left \langle g, h \mid ghg=hgh \right \rangle.
\]
The peripheral subgroup for $\Gamma_M$ is
\[ 
	P_M=\langle g, g^{-4}hg^2h\rangle
\]
with $g$ a standard meridian and $g^{-4}hg^2h$ a standard longitude.

\begin{figure}[h]
    \begin{center}
	\hspace{-1.3cm}
    \subfigure[Generators;]{%
        \label{fig:TrefoilGenerators}%
        \includegraphics[height=3.5cm]{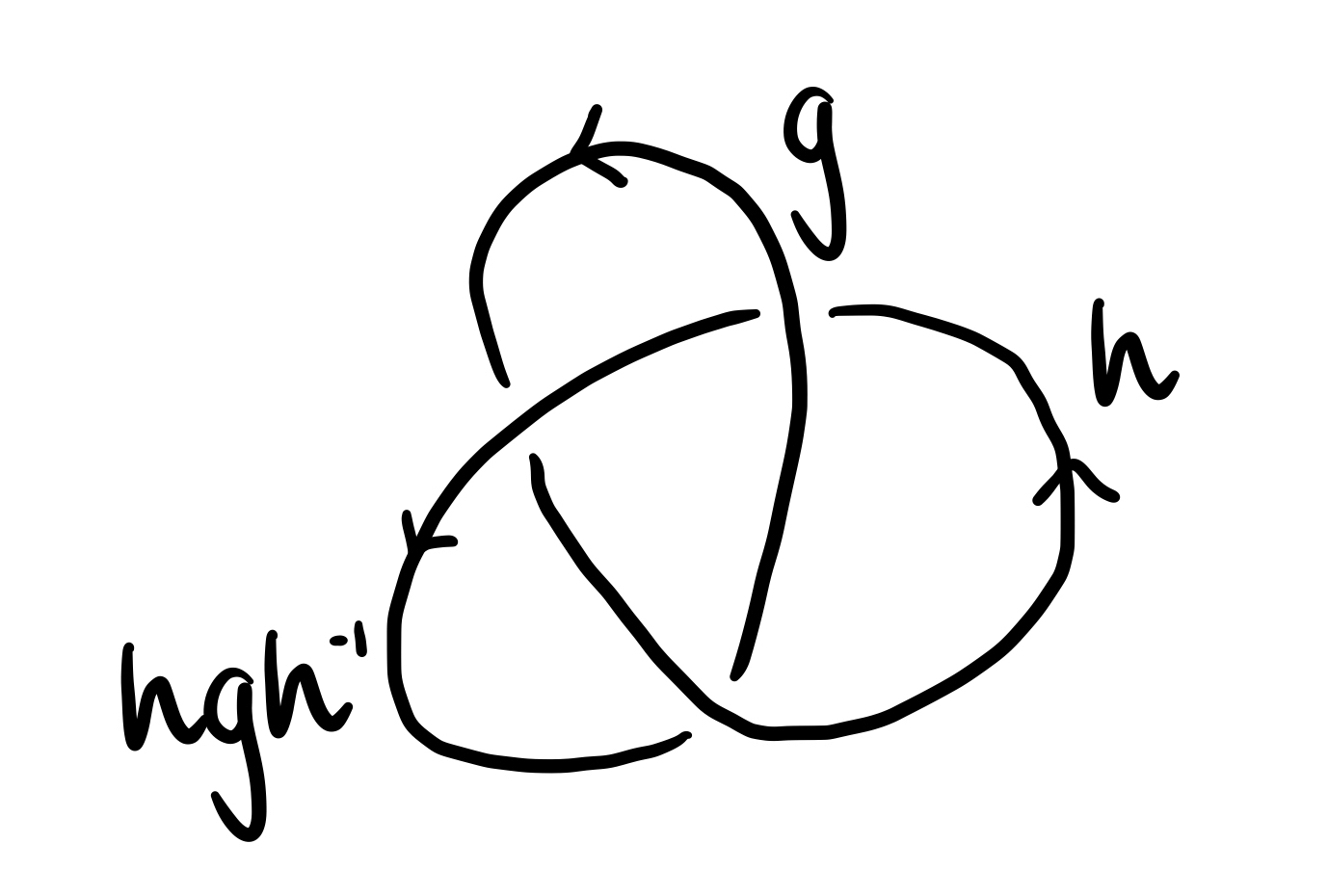}}
 	\quad \quad
        \hspace{0.5cm}
     \subfigure[Seifert surface;]{%
        \label{fig:TrefoilSeifert}%
        \includegraphics[height=3.5cm]{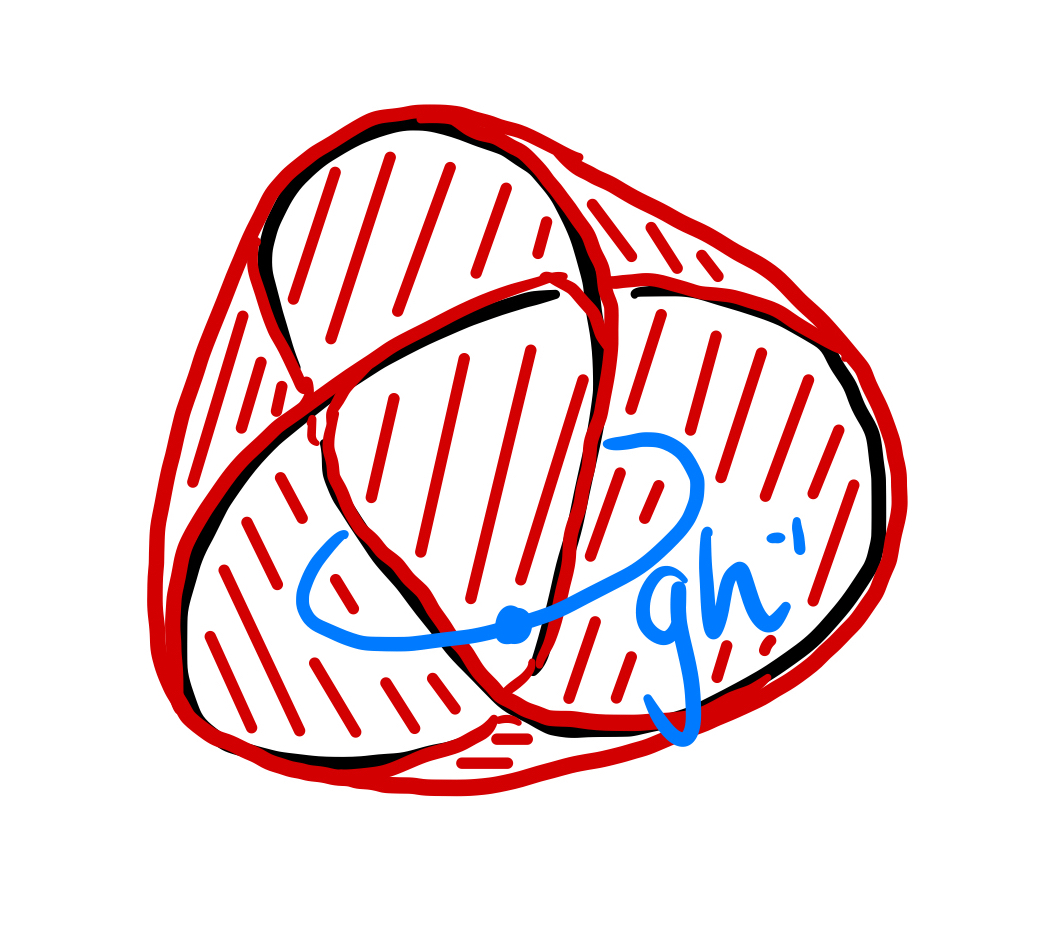}}
        \hspace{1.2cm}
    \\
	\hspace{0.5cm}
 \subfigure[M\"obius band;]{%
        \label{fig:TrefoilMobius}%
        \includegraphics[height=4cm]{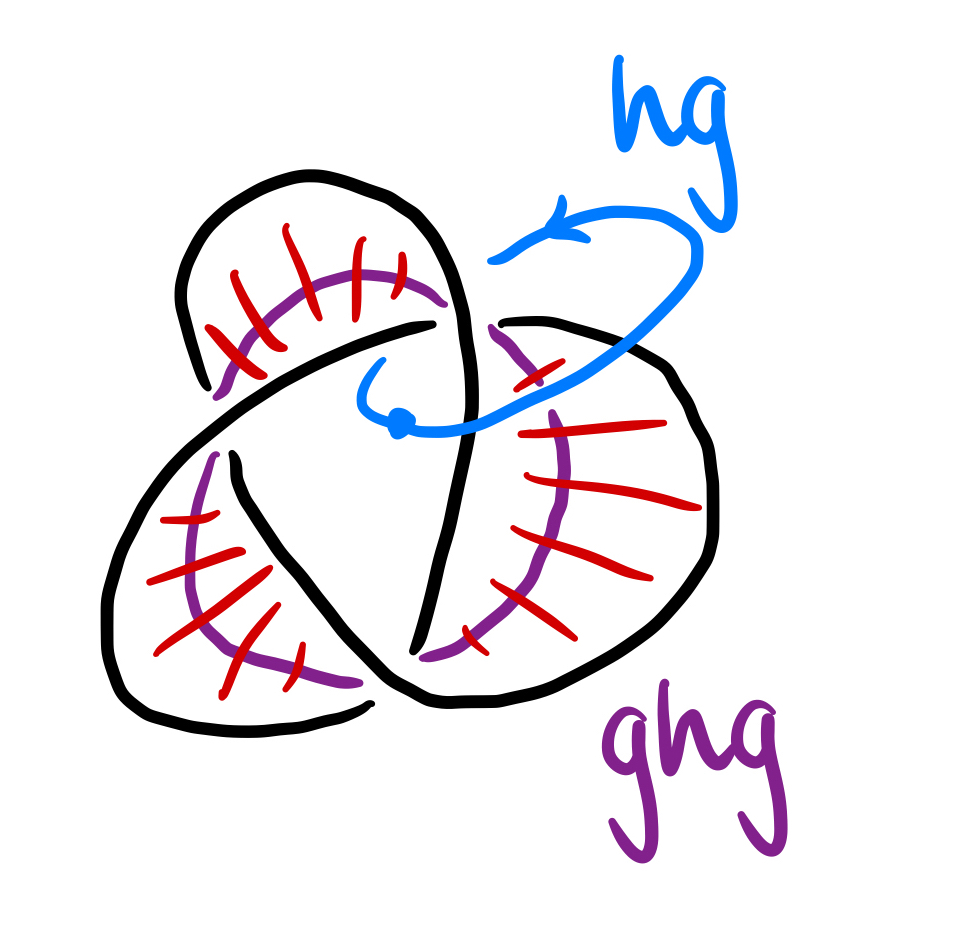}}
        \quad
    \subfigure[Seifert fibration.]{%
        \label{fig:TrefoilD23}%
        \includegraphics[height=4.2cm]{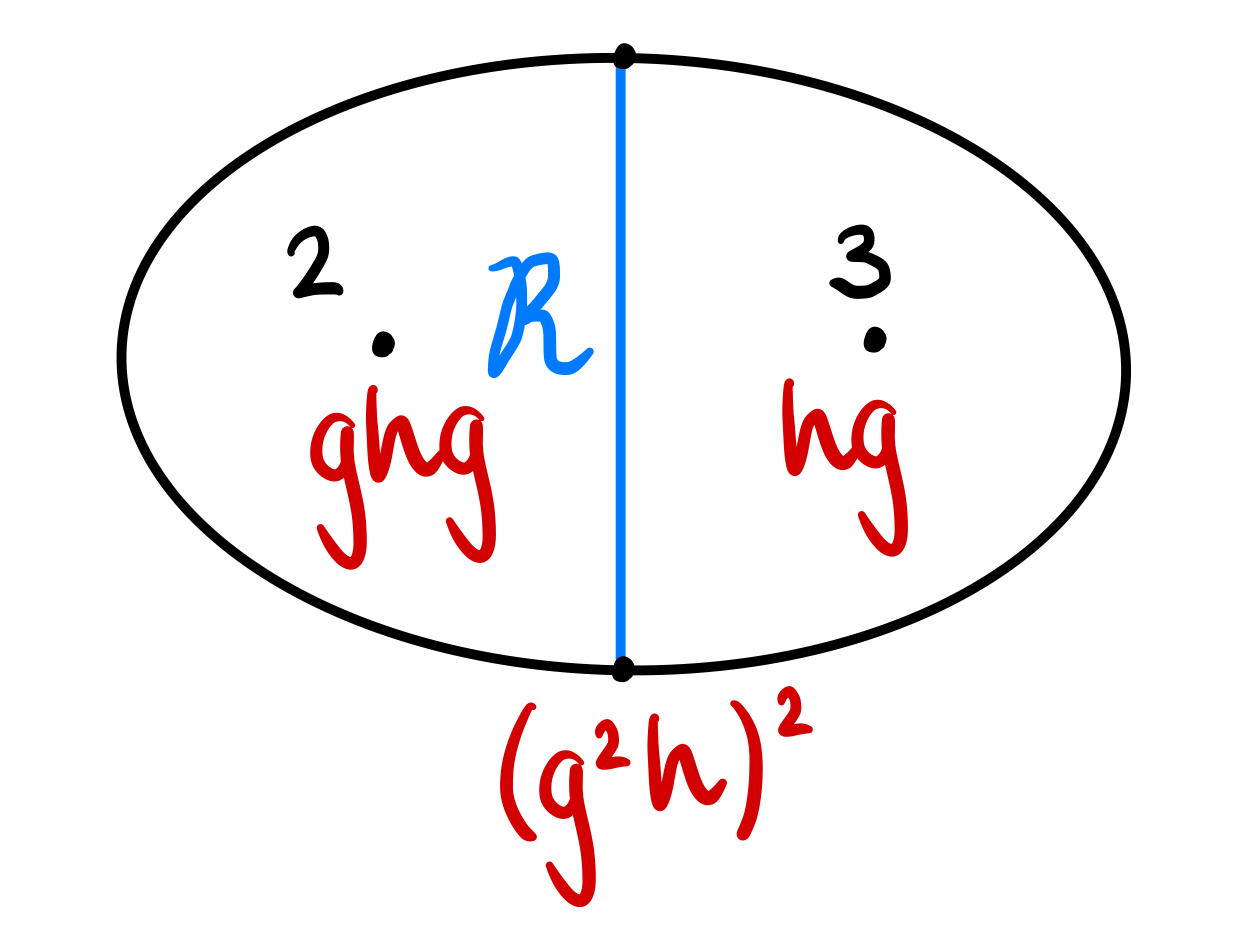}}
        \end{center}
    \caption{Generators of the trefoil complement and elements in complementary regions of the annuli.}
    \label{fig:trefoil}
\end{figure}

As for the previous manifold, we apply \cite[Proposition 1.11]{Hatcher-notes} to show that that there are exactly two connected essential surfaces in $M$. The right-handed trefoil is Seifert fibered with the curve of slope $g^{2}hg^2h$ a regular fibre and base orbifold $D^2(2,3)$. In this case, there is a vertical annulus $\annulusM$ with boundary slope $g^{2}hg^2h$, separating the trefoil knot complement into two solid tori with fibrations $D^2(2)$ and $D^2(3)$ respectively. This is the only connected vertical surface. The unique connected horizontal surface is a once-punctured torus arising from a branched covering of $D^2(2,3).$ This is a Seifert surface $\seifert$ for the trefoil knot with boundary slope $g^{-4}hg^2h.$ Since the slopes of the two surfaces are not parallel, any essential surface in $M$ either consists of parallel copies of $\annulusM$ or of parallel copies of $\seifert.$ Moreover, the surfaces are again uniquely determined by their boundary curves.
  
Consider the variety of characters $\X(M,\F)$. We compute the subvarieties $\X^{\text{irr}}(M,\F)$ and $\X^{\text{red}}(M,\F)$ respectively. 

Each irreducible representation is conjugate to 
\begin{align}
\label{eqn:irredRepnTrefoil2}
	\rho(g) &= \begin{pmatrix} x & 1 \\ 0 & x^{-1}\end{pmatrix}, \quad \rho(h) = \begin{pmatrix} x^{-1} & 0 \\ -1 & x\end{pmatrix}, \quad x\in\F\setminus\{0\}.
\end{align}
By~\Cref{lem:irred_crit}, this representation is irreducible if and only if 
\[
	0 \neq \tr [\rho(g) , \rho(h) ]  - 2 = x^2 - 1 + x^{-2}.
\]
Letting $s = x+x^{-1},$ we map $\rho \mapsto \left(\tr(\rho(g)),\tr(\rho(h)),\tr(\rho(gh))\right)\in \F^3$ giving
\begin{equation*}
	\X^{\text{irr}}(M,\F) = \{(s,s,1)\mid s\in\F\}.
\end{equation*}
This is a curve with a single ideal point. For each $\rho \in \X^{\text{irr}}(M,\F)$, we have $\rho(g^{2}hg^2h) = -E.$ In particular, $\trFunc_{g^{2}hg^2h}$ is constant. Now 
$$\trFunc_{g^{-4}hg^2h} = -(x^6 + x^{-6}) = -s^6 + 6 s^4 - 9 s^2 + 2$$ 
has a pole of order 6 at the ideal point of $\X^{\text{irr}}(M,\F).$ Hence \Cref{pro:detected_surface} implies that $\annulusM$ is detected by this ideal point and $\seifert$ is not. We remark that the order of the pole is related to the fact that the boundary slopes of $\annulusM$ and $\seifert$ have intersection number 6.

The reducible representations are conjugate to the forms
\begin{align}
\label{eqn:redRepnTrefoil0}
	\rho(g) &= \begin{pmatrix} \pm 1 & 1 \\ 0 & \pm 1 \end{pmatrix} = \rho(h),\\
\label{eqn:redRepnTrefoil1}
	\rho(g) &= \begin{pmatrix} x & 0 \\ 0 & x^{-1}\end{pmatrix} = \rho(h), \quad  x\in\F\setminus\{0\},\\
\label{eqn:redRepnTrefoil2} 
	\rho(g) &= \begin{pmatrix} x & 0 \\ 0 & x^{-1}\end{pmatrix}, \quad \rho(h) = \begin{pmatrix} x & 1 \\ 0 & x^{-1}\end{pmatrix}, \quad  x\in\F, \ x^2+x^{-2}=1.
\end{align}
The representations in \Cref{eqn:redRepnTrefoil0,eqn:redRepnTrefoil1} are abelian, and \Cref{eqn:redRepnTrefoil2} defines reducible non-abelian representations.

In characteristic $p=2$, 
\begin{equation*}
\begin{split}
	\X^{\text{red}}(M,\F) 	= \{(s,s,s^2)\mid s\in\F\}.
\end{split}
\end{equation*}
Here, $\X^{\text{irr}}(M,\F)\cap\X^{\text{red}}(M,\F)=\{(1,1,1)\}$. The character in the intersection represents reducible representations conjugate to any of the above forms in \Cref{eqn:irredRepnTrefoil2,eqn:redRepnTrefoil1,eqn:redRepnTrefoil2} with $ x^2+x^{-2}=1$.

In characteristic $p\neq2$,
\begin{equation*}
\begin{split}
	\X^{\text{red}}(M,\F) 	= \{(s,s,s^2-2)\mid s\in\F\}.
\end{split}
\end{equation*}
Here, $\X^{\text{irr}}(M,\F)\cap\X^{\text{red}}(M,\F)=\{(\pm\sqrt{3},\pm\sqrt{3},1)\}$. The characters in the intersection represents reducible representations conjugate to any of the above forms in \Cref{eqn:irredRepnTrefoil2,eqn:redRepnTrefoil1,eqn:redRepnTrefoil2} with $ x^2+x^{-2}=1$.

The Seifert surface $\seifert$ is Poincar\'e dual to the epimorphism $\pi_1(M)\to \Z$ and hence detected by the curve $\X^{\text{red}}(M,\F)$ for any characteristic $p\ge 0.$ This is consistent with the observation that the restriction of $\trFunc_{g^{-4}hg^2h}$ to $\X^{\text{red}}(M,\F)$ is constant equal to 2 (since all characters on $\X^{\text{red}}(M,\F)$ are reducible and the boundary curve of $\seifert$ is a commutator). Moreover, 
$$
\trFunc_{g^2hg^2h} = x^6 + x^{-6} = s^6 - 6 s^4 + 9 s^2 - 2
$$
and hence $\trFunc_{g^2hg^2h}$ has a pole at the ideal point of $\X^{\text{red}}(M,\F).$
Thus \Cref{pro:detected_surface} implies that $\seifert$ is detected by the ideal point of $\X^{\text{red}}(M,\F)$ and $\annulusM$ is not.


\section{Families of graph manifolds and Seifert fibered manifolds}
\label{sec:first_family}

Suppose 
\begin{equation}
	\label{eqn:glueMat}
	\glue=\begin{pmatrix} k & l \\ m & n \end{pmatrix}\in\SL_2(\Z).
\end{equation}
Let $N_\glue$ be the closed orientable 3--manifold obtained by gluing ${\twistKlein}$ to $M$ via $\glue$ so that 
\begin{equation} \label{eqn:gluing}
	g = \left(a^2\right)^k b^l  \quad \text{and} \quad g^{-4}hg^2h = \left(a^2\right)^m b^n.
\end{equation} 
Then $\Gamma_\glue= \pie (N_\glue) = \langle\; a, b, g, h \;\mid\; ghg=hgh, \; aba^{-1}b = 1,\; g = \left(a^2\right)^k b^l, \; g^{-4}hg^2h = \left(a^2\right)^m b^n\; \rangle$.

Note that the abelianisation of $\Gamma_\glue$ is
		\begin{equation*}
		\begin{split}
			H_1(N_\glue)	&= \langle a,b,g,h \mid  ab=ba, b^2=1, g=h, g=a^{2k}b^{l}, 1=a^{2m}b^n \rangle, \\
					&= \langle a,b \mid ab=ba, b^2=1, 1=a^{2m}b^n \rangle,\\
					&= 	\begin{cases} 
								\Z_2\oplus\Z_{2m} & \text{if $n=0\Mod2$}, \\
								\Z_{4m} & \text{if $n=1\Mod2$}
						\end{cases}
		\end{split}
		\end{equation*}
where $\Z_0 = \Z.$ In particular, $H_1(N_\glue)	$ is infinite if and only if $m=0$ and $k=n=\pm1.$


\subsection{Essential surfaces in $N_\glue$}
\label{subsec:surfaces_in_gluing}

We determine all the possible essential surfaces in $N_\glue$. 
\begin{enumerate}
	\item[\surfI] The splitting torus $\splitTorus$ between $M$ and $\twistKlein$. This occurs in $N_\glue$ for all $\glue\in\SL_2(\Z)$ and is an essential separating surface since both $M$ and $\twistKlein$ have incompressible boundary.
\end{enumerate} 

We claim that all other essential surfaces $S$ can be moved via isotopy such that they only have nontrivial intersection with the splitting torus. Indeed, consider $S$ with no intersections with the splitting torus. Then it is a closed, orientable, incompressible surface either in $M$ or in $\twistKlein$. This implies that its components are parallel to $\splitTorus$. Thus every other essential surface $S$ must meet the splitting torus in at least one essential curve. Since we can remove annuli parallel to annuli on $\splitTorus$ by isotopy, it follows that if $S$ is isotoped to have minimal number of intersection curves with $\splitTorus$, then $S\cap (\twistKlein)$ must be an essential surface in $\twistKlein$ and $S\cap M$ must be an essential surface in $M$. The only options for essential surfaces in each manifold $\twistKlein$ and $M$ are given in \Cref{sec:kleinbottle,sec:trefoil} respectively. 

Using this we can list all connected essential surfaces in $N_\glue$ up to isotopy by considering how the boundary curves match up. Since each of $\twistKlein$ and $M$ has two boundary slopes of essential surfaces, this gives four possibilities.
Matching up one boundary curve from each boundary component of $\twistKlein$ and $M$ leaves one degree of freedom in the gluing map. We thus have:

\begin{enumerate}
	\item[\surfII] [$\partial\seifert\leftrightarrow\partial\annulusa]$ 
		Matching up the boundary slope $g^{-4}hg^2h$ with the boundary slope $a^2$ forces $$\glue = \begin{pmatrix} k & \pm1 \\ \mp1 & 0 \end{pmatrix}\in\SL_2(\Z)$$ where $k\in \Z$ is arbitrary. All resulting manifolds $N_\glue$ are graph manifolds and not Seifert fibered. 
		Any copy of $\annulusa$ in $\twistKlein$ matches up with two parallel copies of the Seifert surface in $M$ giving a connected separating genus-2 surface. Any two surfaces obtained this way are isotopic in $N_\glue$ since $M \setminus \seifert \cong \seifert \times (0,1)$. We write $\glueSA=(2\seifert)\cup\annulusa$ for the separating genus-2 surface. 
		
\begin{figure}[t]
    \begin{center}
        \label{fig:glueSA}%
        \includegraphics[height=4.5cm]{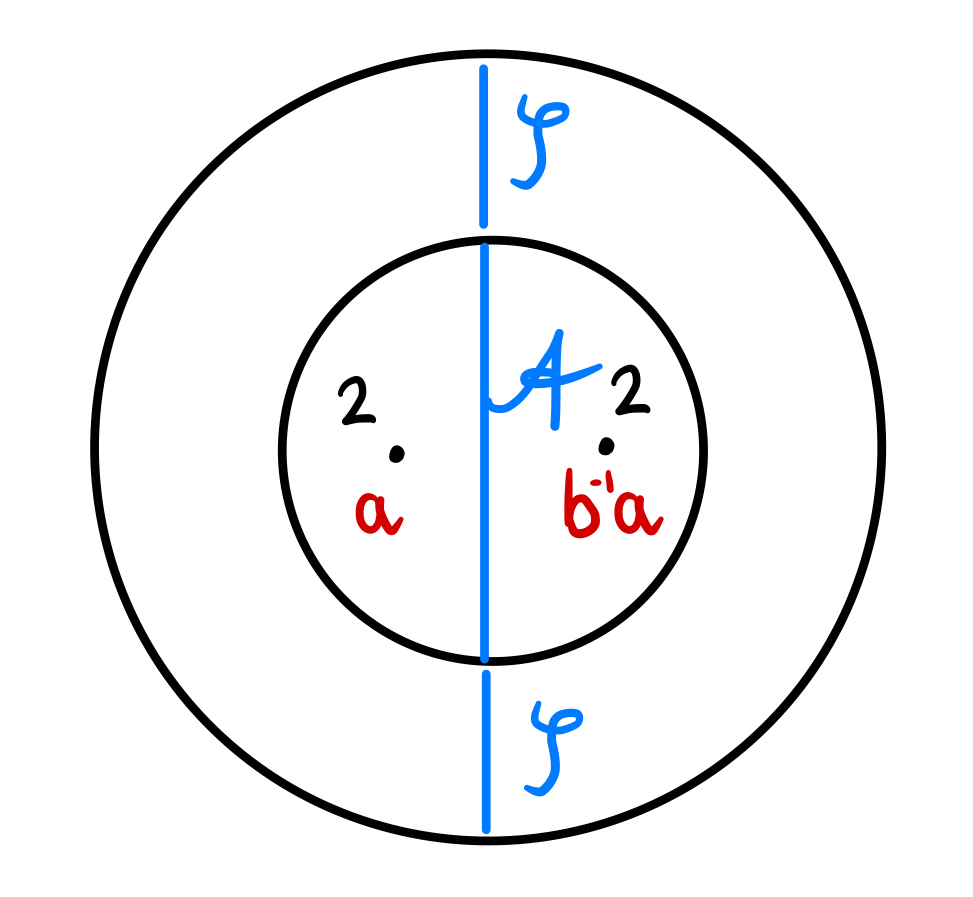}
        \end{center}
    \caption{Separating genus-2 surface arising in $S^2(2,2,2,3)$ from $\glueSA$.}
    \label{fig:SA}
\end{figure}

	\item[\surfIII] [$\partial\seifert\leftrightarrow\partial\annulusb]$ 
		Matching up the boundary slope $g^{-4}hg^2h$ with the boundary slope
 $b$ forces $$\glue = \begin{pmatrix} \pm1 & l \\ 0 & \pm1\end{pmatrix}\in\SL_2(\Z)$$ where $l \in \Z$ arbitrary. All resulting manifolds $N_\glue$ are graph manifolds and not Seifert fibered. 
		Any copy of $\annulusb$ in $\twistKlein$ matches up with two parallel copies of the Seifert surface in $M$ giving a connected surface. Since $\annulusb$ is non-separating in $\twistKlein$, this is a non-separating genus-2 surface, written $\glueSB=(2\seifert)\cup\annulusb$. Any two such surfaces are parallel since $M \setminus \seifert \cong \seifert \times (0,1)$. 
	
	\item[\surfIV] [$\partial\annulusM\leftrightarrow\partial\annulusa]$ 
		Recall that the slopes of $\partial\annulusM$ and $\partial\annulusa$ are regular fibres in Seifert fibrations of $M$ and $\twistKlein$ respectively. The identification forces $$\glue = \begin{pmatrix} k & \pm1 \\ \mp1-6k & \mp6 \end{pmatrix}\in\SL_2(\Z)$$ where $k\in \Z$ is arbitrary. This is a family of Seifert fibered manifolds with base orbifold $S^2(2,2,2,3).$ 
		
		Any combination of surfaces $\annulusM$ and $\annulusa$ will be vertical in the Seifert fibration. Schultens \cite{Schultens-kakimizu-2018} says vertical essential surfaces are in bijective correspondence with isotopy classes of simple closed curves on the four-punctured sphere, which are, in turn, in bijection with isotopy classes of simple closed curves on the torus (see ~\cite[\S2.2.5]{Farb-mapping-2012}). 
		
		Fix meridian $\alpha$ and longitude $\beta$ for the torus. 
		For $q,r\in\Z$, a simple closed curve on the torus is a $(q,r)$-curve if it intersects $\alpha$ $q$ times and $\beta$ $r$ times. We can construct a similar classification of simple closed curves on the four-punctured sphere using the bijection, which proves the simple closed curves are defined by primitive pairs $(q,r)\in\Z^2$. 
		Thus, the vertical essential surfaces are defined by primitive $(q,r)\in\Z^2$ and will correspond to some combination of copies of $\annulusM$ and $\annulusa$. 
		
		We assume that the basis is chosen such that $\splitTorus$ corresponds to the $(1,0)$-curve under the map $N_\glue \to S^2(2,2,2,3)$. Our convention for the $(0,1)$-curve is given in \Cref{fig:RAfamily}.
		We write $\glueRApq=\annulusM\cup\annulusa(q,r)$ for the essential surface $\annulusM\cup\annulusa$ associated with the $(q,r)$-curve. This will be a connected, separating torus and for distinct coprime pairs $(q_1,r_1)$ and $(q_2,r_2)$ the surfaces will be distinct and not isotopic. No horizontal essential surfaces arise in $S^2(2,2,2,3)$.
			
\begin{figure}[t]
    \begin{center}
    \subfigure[$\glueRA{(0,1)}$]{%
        \label{fig:glueRApq0}%
        \includegraphics[height=3.4cm]{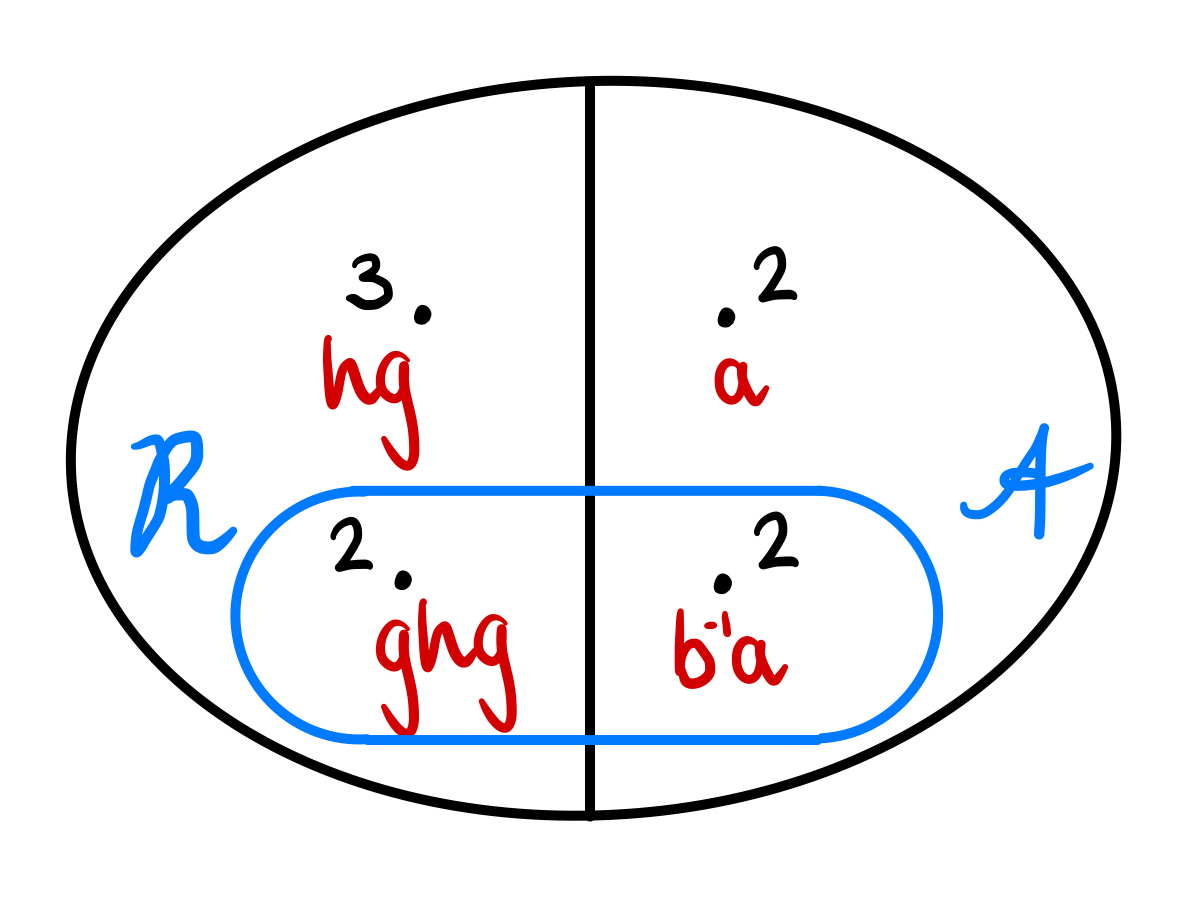}}
     \subfigure[$\glueRA{(1,1)}$]{%
        \label{fig:glueRApq1}%
        \includegraphics[height=3.4cm]{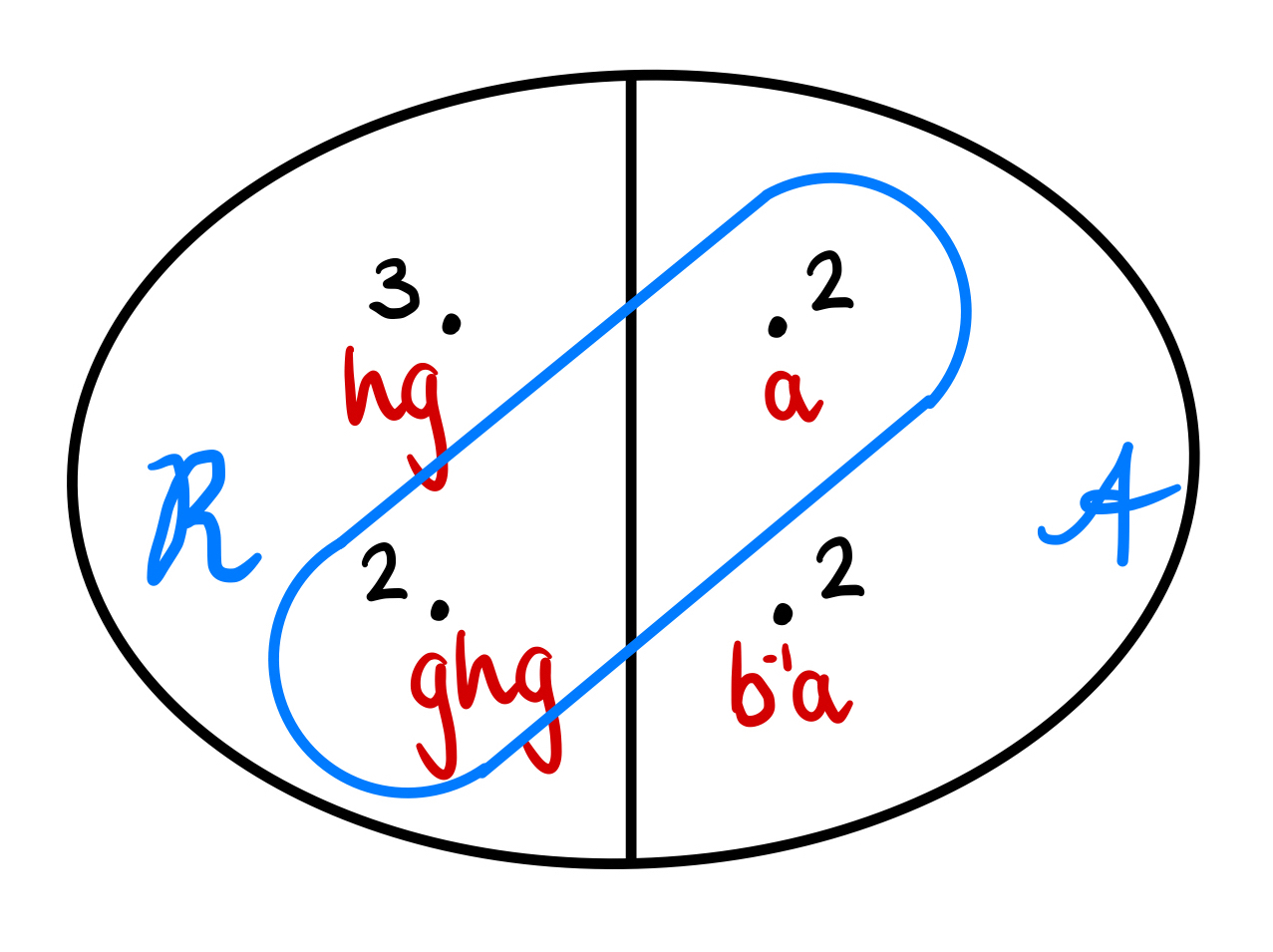}}
     \subfigure[$\glueRA{(2,1)}$]{%
        \label{fig:glueRApq2}%
        \includegraphics[height=3.4cm]{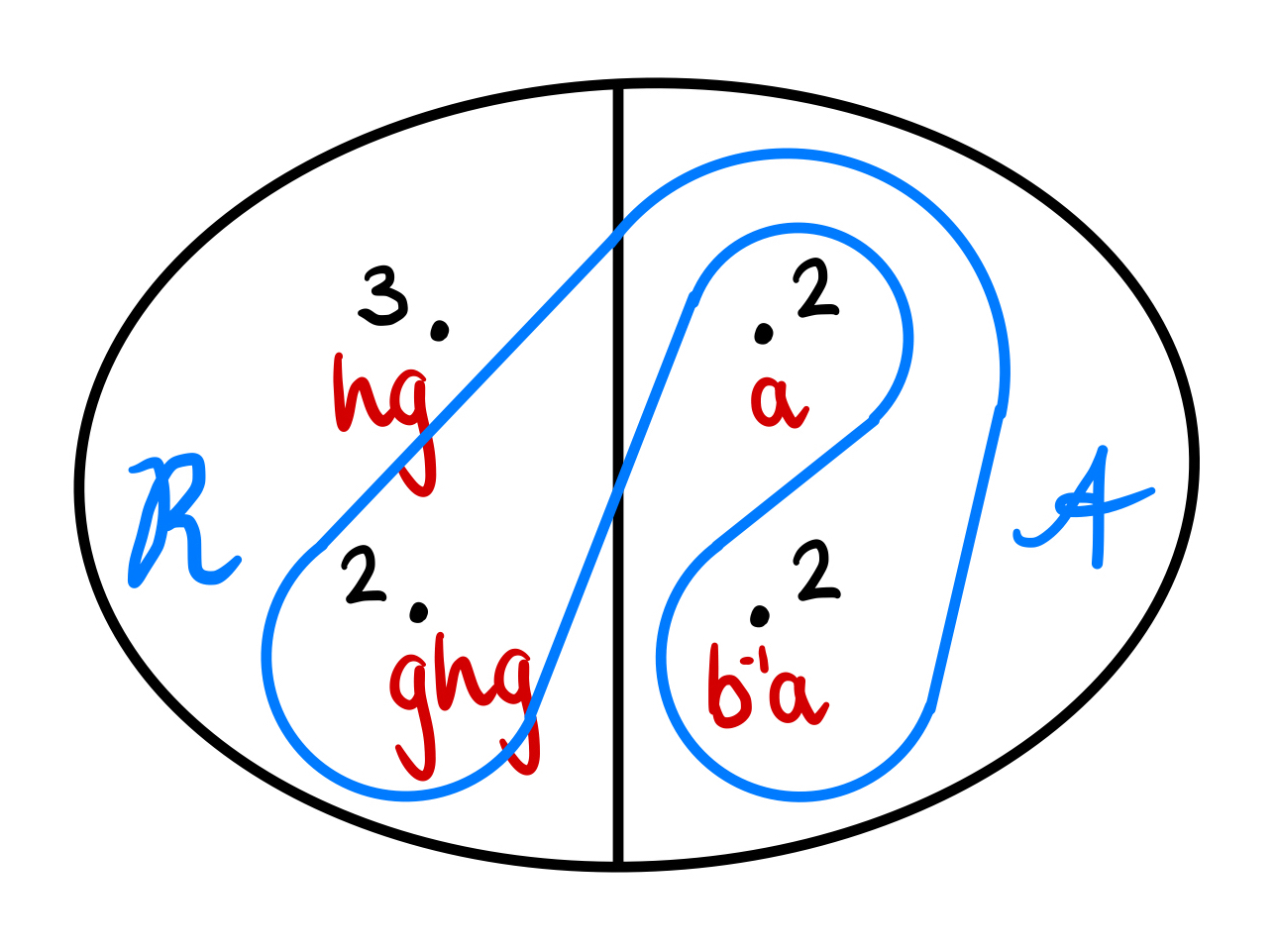}}
        \end{center}
    \caption{Example tori arising in $S^2(2,2,2,3)$ from $\glueRApq$.}
    \label{fig:RAfamily}
\end{figure}

\begin{figure}[t]
    \begin{center}
        \label{fig:glueRApqGeneral}%
        \includegraphics[height=4cm]{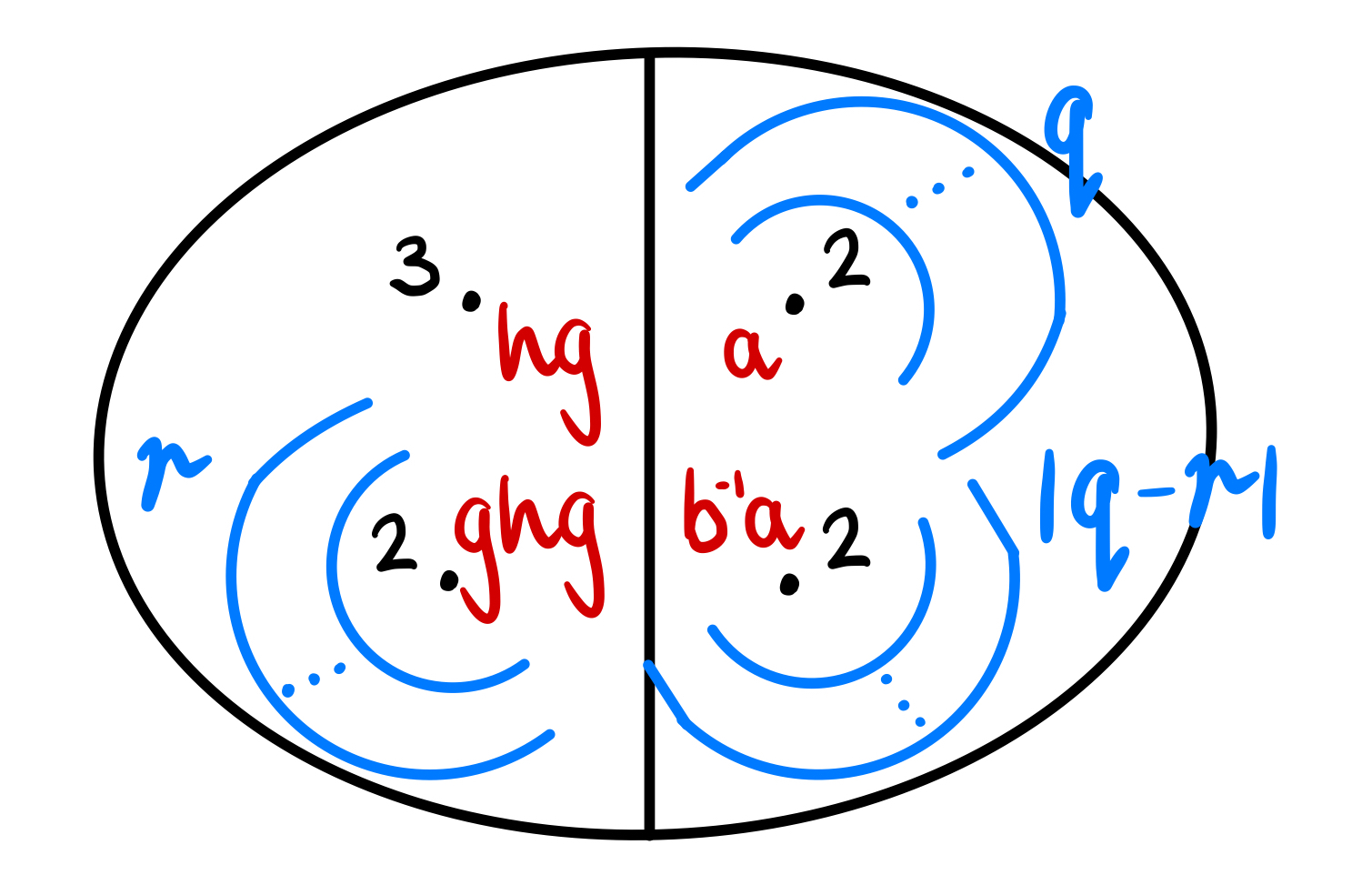}
        \end{center}
    \caption{Example arcs for a general tori arising in $S^2(2,2,2,3)$ from $\glueRApq$. There are $r$ copies of $\annulusM$ and $r$ copies of $\annulusa$ twisted so that $q$ arcs go around the cone point labelled $a$ and $\lvert q-r \rvert$ arcs go around the cone point labelled $b^{-1}a$.}
    \label{fig:RAgeneral}
\end{figure}

	\item[\surfV]  [$\partial\annulusM\leftrightarrow\partial\annulusb]$
	Recall that the slopes of $\partial\annulusM$ and $\partial\annulusb$ are regular fibres in Seifert fibrations of $M$ and $\twistKlein$ respectively. The identification forces
	$$\glue = \begin{pmatrix} \pm1 & l \\ \mp6 & \pm1-6l \end{pmatrix}\in\SL_2(\Z)$$
	where $l\in \Z$ is arbitrary. This is a family of Seifert fibered manifolds with base orbifold $\mathbb{RP}^2(2,3).$ 
	
	Any combination of surfaces $\annulusM$ and $\annulusb$ will be vertical in the Seifert fibration. 
	Note that $\annulusM\cup \annulusb$ gives an embedded Klein bottle in $N_\glue$ and hence is not an essential surface. Taking two copies of each gives a separating torus, written $\glueRB.$ Any two such tori are parallel. Now the core Klein bottle in $\twistKlein$ and the Klein bottle $\annulusM\cup \annulusb$ represent the same element in $H^2(N_\glue, \Z_2).$ However, the tori $\splitTorus$ and $\glueRB=(2\annulusM)\cup (2\annulusb)$ are nevertheless not isotopic. This can be seen by considering the geometric intersection number with a curve represented by $a.$
	No horizontal essential surfaces arise in $\mathbb{RP}^2(2,3)$.
\end{enumerate} 


\subsection{Curves in the variety of characters of $N_\glue$}
\label{subsec:detection_in_glued}

Let $r_\klein \co X(N_\glue,\F) \to X({\twistKlein},\F)$ and $r_M \co X(N_\glue,\F) \to X(M,\F)$ be the restriction maps arising from the inclusions $\pi_1({\twistKlein}) \to \pi_1(N_\glue)$ and $\pi_1(M) \to \pi_1(N_\glue)$ respectively. Suppose $C \subset X(N_\glue,\F)$ is a curve. We have the following possibilities for $C$:

\begin{enumerate}
\item[{\curveI}] Both $r_\klein(C)$ and $r_M(C)$ are curves. It follows from the description of the curves in $X({\twistKlein},\F)$ and $X(M,\F)$ that for each ideal point $\xi$ of $C$, there is a unique primitive class $\alpha \in \im (\pi_1(\splitTorus)\to\pi_1(N_\glue))$ with $v_\xi(\trFunc_\alpha)\ge 0$ and all other primitive classes $\gamma\in \im (\pi_1(\splitTorus)\to\pi_1(N_\glue))$ satisfy $v_\xi(\trFunc_\gamma)< 0.$ Hence none of these classes $\gamma$ are contained in an edge stabiliser of the action on the Bass-Serre tree associated with the ideal point, and hence none of these are homotopic to a curve on the surface detected by $\xi.$ It follows that any essential surface detected by $C$ is as listed in \surfII--\surfV. Except for the case \surfIV, the detected surface is uniquely determined by this information (up to taking parallel copies). 
\item[{\curveII}]  One of $r_\klein(C)$ and $r_M(C)$ is a curve and the other is a point. In this case, the curve again forces that for each  ideal point $\xi$ of $C$ there is a primitive class $\gamma\in \im (\pi_1(\splitTorus)\to\pi_1(N_\glue))$ with $v_\xi(\trFunc_\gamma)< 0.$ However, since the other projection is a point, all representations with character on $C$ have constant trace functions $\trFunc_\beta$ for all $\beta\in \im (\pi_1(\splitTorus)\to\pi_1(N_\glue)).$ This implies $v_\xi(\trFunc_\gamma)\ge 0,$ a contradiction. 
\item[{\curveIII}]  Each of $r_\klein(C)$ and $r_M(C)$ is a point. In this case, all representations with character on $C$ have constant trace functions $\trFunc_\gamma$ for all $\gamma\in \im (\pi_1(\splitTorus)\to\pi_1(N_\glue)),$ for all $\gamma\in \im (\pi_1(\twistKlein)\to\pi_1(N_\glue)),$ and for all $\gamma\in \im (\pi_1(M)\to\pi_1(N_\glue)).$ It follows that we may choose the splitting torus $\splitTorus$ to be an essential surface detected by $\xi.$ To rule out that any of the essential surfaces $S$ listed in \surfII--\surfV\  is also detected by $C$, according to \Cref{pro:stabilisers} we need to exhibit an element $\gamma \in \im (\pi_1(N'_\glue)\to\pi_1(N_\glue))$ with $v_\xi(\trFunc_\gamma)< 0$ for $N'_\glue$ a component of $N_\glue\setminus S.$
\end{enumerate}

In the following sections we determine when essential surfaces are detected by curves in the variety of characters of $N_\glue$, broken into cases for detecting: the splitting torus as in \surfI; additional surfaces in the graph manifold as in \surfII\ and \surfIII; additional surfaces in the Seifert fibration with base orbifold $S^2(2,2,2,3)$ as in \surfIV; and additional surfaces in the Seifert fibration with base orbifold $\mathbb{RP}^2(2,3)$ as in \surfV.

The approach for proving each result is similar.
For ease of notation write 
$$\left(A,B,G,H\right)=\left(\rho(a),\rho(b),\rho(g),\rho(h)\right)\in\SLF^4.$$ 

We analyse the variety of characters of $N_\glue$ for each case and outline when a curve appears, up to conjugacy. 
We start with the options for representations of ${\twistKlein}$ outlined in \Cref{sec:kleinbottle}. To extend the representations to $N_\glue$, the matrices must satisfy the gluing equations
\begin{equation}
\label{eqn:gluingMat}
	G=A^{2k}B^{l},\quad G^{-4}HG^2H=A^{2m}B^{n}.
\end{equation}

If $B=\pm\idMat$, as in \Cref{eqn:redRepnKlein1}, the equations simplify to 
\begin{equation*}
	G=A^{2k}(\pm1)^l,\quad HA^{4k}H=(\pm1)^nA^{8k+2m}.
\end{equation*}

If $B\neq\pm \idMat$, from \Cref{prop:kleinbottle} either $\rho$ factors through the abelianisation of $N_\glue$ or $A^2=\pm \idMat$. We have $A^2=\pm \idMat$ in the case of \Cref{eqn:irredRepnKlein,eqn:redRepnKlein2}. Here, the equations simplify to
\begin{equation*}
	G=(\pm1)^kB^{l},\quad HB^{2l}H=(\pm1)^m B^{4l+n}.
\end{equation*}

The resulting matrices must also satisfy group relation for $M$,
\begin{equation}
\label{eqn:trefoilRelationMat}
	GHG=HGH.
\end{equation}

The general procedure is as follows. Matrix $G$ is defined by \Cref{eqn:gluingMat}. Write 
\begin{equation}
\label{eqn:h}
	H=\begin{pmatrix} e & f \\ g & h \end{pmatrix}, \quad eh-fg=1.
\end{equation}

Expanding the powers of $A$ and $B$, \Cref{eqn:gluingMat} gives four equations for the entries of $H$ in terms of the other variables in addition to the determinant condition. 
Outside of the cases where there is a curve, we either find a contradiction in these equations, a contradiction when considering the group relation \Cref{eqn:trefoilRelationMat}, a contradiction to the gluing matrix $\glue\in\SL_2(\Z)$, or that all variables are specified and so we have at most some 0--dimensional components.

Where a curve is possible, we write for ease of notation $t_\gamma=I_\gamma(\rho)=\tr (\rho(\gamma))$ for $\gamma \in \pi_1(N_\glue)$ and $\rho\in\R(N_\glue,\F)$. The variety $\X_p(N_\glue)$ has as coordinates the traces of ordered single, double and triple products of the generators in any characteristic (see~\Cref{cor:trace_function_characterisation}).
These coordinates are
\[
	\left(t_a,t_b,t_g,t_h,t_{ab},t_{ag},t_{ah},t_{bg},t_{bh},t_{gh},t_{abg},t_{abh},t_{agh},t_{bgh}\right) \in \F^{14}.
\]

In order to determine what surface is detected by a curve, we appeal to \Cref{pro:stabilisers,pro:detected_surface}, the discussions of essential surfaces in $\twistKlein$ and $M$ in \Cref{sec:kleinbottle,sec:trefoil}, the classification of connected essential surfaces in $N_\glue$ given by \surfI--\surfV, and the classification of the types of curves that appear in \curveI--\curveIII.


\subsection{Detecting the splitting torus}
\label{subsec:detection_torus}

The manifold $N_\glue$ always contains the splitting torus $\splitTorus$ for all $\glue \in\SL_2(\Z)$. We get two results for when the splitting torus is detected over characteristic $p\neq2$ and over characteristic 2.
\begin{lemma} 
\label{lem:gluingCurveTorus1} 
	Let $\F$ be an algebraically closed field of characteristic $p\neq 2$. 
	There is a curve in $\X(N_\glue, \F)$ that detects the splitting torus $\splitTorus$ if and only if $\glue$$=$$\begin{pmatrix} k & \pm1 \\ \mp1-6k & \mp6 \end{pmatrix} \in\SL_2(\Z)$ for $k\in\Z$. Moreover, there are infinitely many such curves whenever $\glue$ is of this form. 
\end{lemma}

\begin{lemma} 
\label{lem:gluingCurveTorus2} 
	Let $\F$ be an algebraically closed field of characteristic $2$. There is a curve in $\X(N_\glue, \F)$ that detects the splitting torus $\splitTorus$ if and only if $\glue = \begin{pmatrix} k & 1 \\ 1 & 0 \end{pmatrix} \Mod2$ $\in\SL_2(\Z\slash2\Z)$ for $k\in\{0,1\}$.

Moreover, the curve is unique unless $\glue$$=$$\begin{pmatrix} k & \pm1 \\ \mp1-6k & \mp6 \end{pmatrix} \in\SL_2(\Z)$ for $k\in\Z$, in which case there are infinitely many curves.
\end{lemma}

\Cref{thm:graph-torus-never} now follows by taking (for instance) the family 
$\glue_k = \begin{pmatrix} 1+6k & k \\ 6 & 1 \end{pmatrix}$ and 
\Cref{thm:graph-torus-char2} by taking (for instance) the family
$\glue'_k = \begin{pmatrix} k & 1+6k \\ -1 & -6 \end{pmatrix}$ for $k \in \mathbb{N}.$
These families have been chosen such that no other essential surfaces are present. We have $\glue_k$ clearly does not satisfy the requirement for \Cref{lem:gluingCurveTorus1,lem:gluingCurveTorus2} and $\glue'_k$ only satisfies \Cref{lem:gluingCurveTorus2} but not \Cref{lem:gluingCurveTorus1}.

\begin{proof}[Proof of \Cref{lem:gluingCurveTorus1}]
Assume $p\neq2$. We follow the approach described in \Cref{subsec:detection_in_glued}.
The splitting torus $\splitTorus$ can only be detected by a curve of type \curveIII\  where the restrictions to the variety of characters of $\twistKlein$ and $M$ are points.  
We will see that there are two ways to get such a curve in characteristic $p\neq 2$, one of which comes from a 2--dimensional component that will give infinitely many curves. 

Consider the options for matrix pairs $\left(A,B\right)$ from \Cref{eqn:irredRepnKlein,eqn:redRepnKlein1,eqn:redRepnKleinCharNon22} in turn. 

We find one way to extend the representation from \Cref{eqn:irredRepnKlein} to $N_\glue$ with a curve of type \curveIII.
For gluing matrix $\glue=\begin{pmatrix} k & \pm1 \\ \mp1-6k & \mp6 \end{pmatrix}\in\SL_2(\Z)$ we get a 2--dimensional component of representations
\begin{equation*}
\begin{split}
	A = \begin{pmatrix} 0 & -1 \\ 1 & 0\end{pmatrix}, 
	&\quad 
	B = \begin{pmatrix} y & 0 \\ 0 & y^{-1}\end{pmatrix}, \\ 
	G = (-1)^k\begin{pmatrix} y^{\pm1} & 0 \\ 0 & y^{\mp1}\end{pmatrix}, 
	&\quad 
	H = \begin{pmatrix}\frac{\mp(-1)^{k}y^{\mp2}}{y-y^{-1}} & f \\  \frac{-y^2+1-y^{-2} }{f(y-y^{-1})^2}  & \frac{\pm(-1)^{k} y^{\pm2}}{y-y^{-1}} \end{pmatrix}, \quad f,y\in\F\setminus\{0\}.
\end{split}
\end{equation*}

The corresponding traces give a 2-dimensional component of the variety of characters, 
\begin{equation}
\begin{split}
\label{eqn:cTorusI}
	\cTorusI =	&\left\{\left(0, y+y^{-1}, (-1)^k(y+y^{-1}), (-1)^k (y+y^{-1}), 0, 0, \frac{y^2 - 1 + y^{-2}}{f (y-y^{-1})^2}+f,    \right.\right.\\
			&\hspace{1cm}\left.\left.(-1)^k(y^{1\pm1}+y^{-1\mp1}),  \pm(-1)^k\frac{y^{-1\pm2}-y^{1\mp2}}{y-y^{-1}},{1},0, \frac{y^2 - 1 + y^{-2}}{f y(y-y^{-1})^2}+fy, \right. \right. \\
			&\hspace{1cm}\left.\left. (-1)^k\left(\frac{y^2 - 1 + y^{-2}}{fy^{\pm1}(y-y^{-1})^2}+fy^{\pm1}\right), \pm\frac{y^{-1\pm1}-y^{1\mp1}}{y-y^{-1}}\right)\ \Big|\ f,y\in\F\setminus\{0\}\right\}
\end{split}
\end{equation}
for $l=\pm1$ and $k\in\Z$ even or odd. 
		
Fix $y\in\F\setminus\{0\}$ to be constant. This defines curves $C_y^{\pm}(k)$ inside $\cTorusI$ since $I_{ah}$ is not constant. The restrictions to $\twistKlein$ and $M$ are points, 
\begin{align*}
	r_\klein\left(C_y^{\pm}(k)\right) 	&= (0,y+y^{-1},0) \in \X^{\text{irr}}(\twistKlein,\F_p),\\
	r_M\left(C_y^{\pm}(k)\right) 	&= \left\{ \left((-1)^k(y+y^{-1}),(-1)^k(y+y^{-1}),1\right)\mid k\in\Z \right\} \subset \X^{\text{irr}}(M,\F_p),
\end{align*}
so it is of is of type \curveIII.
There are infinitely many such curves for $y\in\F\setminus\{0\}$ that are distinct for each $y+y^{-1}$.
We deduce the curve detects $\splitTorus$ whenever $\glue=\begin{pmatrix} k & \pm1 \\ \mp1-6k & \mp6 \end{pmatrix}\in\SL_2(\Z)$.

We find no ways to extend the representation from \Cref{eqn:redRepnKlein1} to get a curve of type \curveIII\ when $p\neq2$. 

We find one way to extend the representation from \Cref{eqn:redRepnKleinCharNon22} to $N_\glue$ with a curve of type \curveIII.
For gluing matrix {$\glue=\begin{pmatrix} k & \pm1 \\ \mp1-6k & \mp6 \end{pmatrix}\in\SL_2(\Z)$} we get the curves of representations
\begin{equation*}
\begin{split}
	A = \begin{pmatrix} x & 0 \\ 0 & -x\end{pmatrix}, &\quad 
	B = \begin{pmatrix} s & 1 \\ 0 & s \end{pmatrix}, \quad x\in\F\setminus\{0\},\; x^2=-1 \\ 
	G = (-1)^k\begin{pmatrix} s & \pm1 \\ 0 & s \end{pmatrix}, &\quad 
	H = \begin{pmatrix} e & \mp(-1)^k(-1 + 2 s (-1)^{k} e - e^2) \\ \mp(-1)^k & 2 s (-1)^{k}- e\end{pmatrix},\; s\in\{-1,1\}, e\in\F. 	
\end{split}
\end{equation*}
The condition on $\glue$ arises again from the gluing equation, which requires $n=-6l$ when $p\neq2$. 
			
The corresponding traces give curves 		
\begin{equation}
\begin{split}
\label{eqn:cTorusIII}
	{\cTorusIII} 	&=\left\{ \left(0, 2 s, 2 (-1)^k s, 2 (-1)^k s, 0, 0, \mp(-1)^{k} + 2 (e - (-1)^{k} s) x, 2 (-1)^k, (-1)^{k} (2 \mp1), 1, 0, \right. \right. \\
				&\hspace{0.8cm} \left. \left. 2 (-(-1)^{k} + e s) x \mp (-1)^{k} (s + x), -(-1)^{k} (\pm(-1)^k s + 3 (-1)^k x - 2 e s x), s(1\mp1)\right) \ \Big|\ e\in \F\right\}
\end{split}
\end{equation}
for $k\in\Z$ even or odd, $l=\pm1$, $s=\pm1$, and $x\in\F\setminus\{0\},\ x^2=-1$.

These clearly define curves since $I_{ah}$ is not constant. The restrictions to $\twistKlein$ and $M$ are points, 
\begin{align*}		
	r_\klein\left(\cTorusIII\right)	&= \left\{(0,2s,0)\ | \ s=\pm1 \right\}\subset X^\text{irr}({\twistKlein},\F_p)\cap X^\text{red}({\twistKlein},\F_p),\\
	r_M\left(\cTorusIII\right) 		&= \left\{(2(-1)^ks,2(-1)^ks,1)\ | \ k\in\Z,\ s=\pm1 \right\}\in X^\text{irr}(M,\F_p),
\end{align*}
so it is of type \curveIII.
We deduce the curves $\cTorusIII$ detect $\splitTorus$ whenever $\glue= \begin{pmatrix} k & \pm1 \\ \mp1-6k & \mp6 \end{pmatrix}\in\SL_2(\Z)$.

No other curves of type \curveIII\ are found for any gluing matrix and other representations when $p\neq2$. This proves the result.
\end{proof}

The proof of the result over characteristic 2 is similar. 
\begin{proof}[Proof of \Cref{lem:gluingCurveTorus2}]
Assume $p=2$.
We follow the approach described in \Cref{subsec:detection_in_glued}.
The splitting torus $\splitTorus$ can only be detected by a curve of type \curveIII\  where the restrictions to the variety of characters of $\twistKlein$ and $M$ are points.  
We will see that there are two ways to get such a curve in characteristic 2, 
{one arises when $\glue=\begin{pmatrix} k & \pm1 \\ \mp1-6k & \mp6 \end{pmatrix}\in\SL_2(\Z)$, which comes from a 2--dimensional component that will give infinitely many curves, and the other arises when $\glue =  \begin{pmatrix} k & 1 \\ 1 & 0 \end{pmatrix} \Mod2 \in\SL_2(\Z\slash2\Z)$.}

Consider the options for matrix pairs $\left(A,B\right)$ from \Cref{eqn:irredRepnKlein,eqn:redRepnKlein1,eqn:redRepnKleinChar21} in turn. 

We find one way to extend the representation from \Cref{eqn:irredRepnKlein} to $N_\glue$ with a curve of type \curveIII\ using the 2--dimensional component $\cTorusI$ given in \Cref{eqn:cTorusI} found in the proof of \Cref{lem:gluingCurveTorus1}. The same result follows and we get there are infinitely many such curves that detect $\splitTorus$ whenever $\glue=\begin{pmatrix} k & \pm1 \\ \mp1-6k & \mp6 \end{pmatrix}\in\SL_2(\Z)$.

We find no ways to extend the representation from \Cref{eqn:redRepnKlein1} to get a curve of type \curveIII\ when $p=2$. 

We find one way to extend the representation from \Cref{eqn:redRepnKleinChar21} to $N_\glue$ with a curve of type \curveIII.
For gluing matrix is $\glue=\begin{pmatrix} k & l \\ m & 2q \end{pmatrix}\in\SL_2(\Z)$, we get the curve of representations
\begin{equation*}
\begin{split}
	A = \begin{pmatrix} 1 & u \\ 0 & 1\end{pmatrix}, &\quad 
	B = \begin{pmatrix} 1 & 1 \\ 0 & 1\end{pmatrix}, \\
	G = \begin{pmatrix} 1 & 1 \\ 0 & 1 \end{pmatrix}, &\quad 
	H = \begin{pmatrix} 1 & 0 \\ 1 & 1\end{pmatrix}, \quad u\in\F.
\end{split}
\end{equation*}

Note the determinant condition in $\Z$ forces $l=1\Mod2$ and $m=1\Mod2$. 
The condition on $\glue$ arises from the gluing equation. We have 
\[
	HB^{2l}H-B^{n+4l}=\begin{pmatrix} e(e+h) & f(e+h) +n \\ g(e+h) & h(e+h)\end{pmatrix}
\]
evaluating to the zero matrix, which forces $n=0\Mod2$.
			
The representation is conjugate to the form, 	
\begin{equation}
\label{eqn:char2iii}
\begin{split}
	A = \begin{pmatrix} 1 & u \\ 0 & 1\end{pmatrix}, &\quad 
	B = \begin{pmatrix} 1 & 1 \\ 0 & 1\end{pmatrix}, \\ 
	G = \begin{pmatrix} 1 & 1 \\ 0 & 1 \end{pmatrix}, &\quad 
	H = \begin{pmatrix} e & 1+e^2 \\ 1 & e \end{pmatrix}, \quad u\in\F,
\end{split}
\end{equation}
which is useful when comparing to the $p\neq 2$ case.
	
The corresponding traces give the curve
\begin{equation}
\label{eqn:cTorusII}
	\cTorusII=\left\{\left(0,0,0,0,0,0,u,0,1,1,0,u+1,u+1,0\right)\mid u\in\F \right\}.
\end{equation}	
This clearly defines a curve since $I_{ah}$ is not constant. The restrictions to $\twistKlein$ and $M$ are points, 
\begin{align*}
	r_\klein\left(\cTorusII\right) 	&= (0,0,0)\in X^\text{irr}({\twistKlein},\F_2)\cap X^\text{red}({\twistKlein},\F_2),\\
	r_M\left(\cTorusII\right) 		&= (0,0,1)\in X^\text{irr}(M,\F_2),
\end{align*}
so it is of is of type \curveIII.
We deduce the curve detects $\splitTorus$ whenever $\glue = \begin{pmatrix} k & 1 \\ 1 & 0 \end{pmatrix} \Mod2$ $\in\SL_2(\Z\slash2\Z)$ for $k\in\{0,1\}$.

The representations and associated curves correspond to $\cTorusIII$ in \Cref{eqn:cTorusIII} found when $p\neq2$ in the proof of \Cref{lem:gluingCurveTorus1}, just for more restrictive gluing matrices. 
Reducing the above representations modulo 2 is the same as the representation in \Cref{eqn:char2iii} with $u=0$ and reducing the curve $\cTorusIII$ modulo 2 is the same as the curve $\cTorusII$ with $u=1$. 
We find that if $p\neq2$ the remaining degree of freedom, $e$, cannot be removed by conjugation and is in fact necessary to find a curve in the variety of characters.

No other curves of type \curveIII\ are found for any gluing matrix and other representations when $p=2$. This proves the result.
\end{proof}


\subsection{Curves in the graph manifolds}
\label{subsec:detection_in_graph}

The manifold $N_\glue$ forms a graph manifold if and only if $$\glue = \begin{pmatrix} k & l \\ m & n\end{pmatrix} \in\SL_2(\Z),\quad m\neq\pm6\neq n.$$ 
The classification of essential surfaces in $N_\glue$ in \Cref{subsec:surfaces_in_gluing} shows only the cases \surfI, \surfII, and \surfIII\ can appear. The case \surfI\ is already covered in \Cref{lem:gluingCurveTorus1,lem:gluingCurveTorus2}. The other cases \surfII\ and \surfIII\ and related curves in the variety of characters of the graph manifold are covered in the next result. 

\begin{lemma} 
\label{lem:gluingCurveGraph} 
Let $\F$ be an algebraically closed field of characteristic $p$. Consider $N_\glue$ a graph manifold with gluing matrix $\glue=\begin{pmatrix} k & l \\ m & n \end{pmatrix}\in\SL_2(\Z)$, $m\neq\pm6\neq n$. Curves that detect \surfII\ and \surfIII\ in $\X(N_\glue,\F)$ are characterised by the following.  
\begin{itemize}
	\item In case \surfII, for $\glue = \begin{pmatrix} k & \pm1 \\ \mp1 & 0 \end{pmatrix}\in\SL_2(\Z)$, there is a curve that detects $\glueSA$ if $p=2$ and no curve that detects $\glueSA$ if $p\neq 2$;
	\item In case \surfIII, for $\glue = \begin{pmatrix} \pm1 & l \\ 0 & \pm1 \end{pmatrix}\in\SL_2(\Z)$, there is a curve that detects $\glueSB$.
\end{itemize}

In particular, $\glueSA$ is only detected in characteristic $2$ and $\glueSB$ is detected whenever it appears. 
\end{lemma}

Recall from \Cref{lem:gluingCurveTorus1,lem:gluingCurveTorus2}, the splitting torus $\splitTorus$ is detected in case \surfII\ in characteristic $p=2$ and is never detected otherwise. Thus, \Cref{lem:gluingCurveTorus1,lem:gluingCurveTorus2,lem:gluingCurveGraph}  give the statement of \Cref{thm:graph-torus-gen2}.
Similarly, \Cref{lem:gluingCurveTorus1,lem:gluingCurveTorus2,lem:gluingCurveGraph}  give the statement of \Cref{thm:graph-torus-nonsepgen2}.

\begin{proof} 
We follow the approach described in \Cref{subsec:detection_in_glued}. The essential surface $\glueSA$ can only be detected by either a curve of type \curveI\ where the restriction to the variety of characters of $\twistKlein$ detects $\seifert$ and the restriction to the variety of characters of $M$ detects $\annulusa$ or a curve of type \curveIII\ that has non-negative valuation on all simple closed curves in components of $N_\glue\setminus\glueSA$.
Similarly, the essential surface $\glueSB$ can only be detected by either a curve of type \curveI\ where the restriction to the variety of characters of $\twistKlein$ detects $\seifert$ and the restriction to the variety of characters of $M$ detects $\annulusb$ or a curve of type \curveIII\ that has non-negative valuation on all simple closed curves in components of $N_\glue\setminus\glueSB$.

Consider first the options for matrix pairs $\left(A,B\right)$ from \Cref{eqn:irredRepnKlein,eqn:redRepnKlein1,eqn:redRepnKleinChar21,eqn:redRepnKleinCharNon22} in turn that could give appropriate curves of type \curveI\ that detect $\glueSA$. 

For characteristic $2$ we find one way to extend the representation from \Cref{eqn:irredRepnKlein} to $N_\glue$ with a curve of type \curveI.
For gluing matrix $\glue=\begin{pmatrix} k & \pm1 \\ \mp1 & 0 \end{pmatrix}\in\SL_2(\Z)$ we get the curve of representations
\begin{equation*}
\begin{split}
	A = \begin{pmatrix} 0 & 1 \\ 1 & 0\end{pmatrix}, 
	&\quad 
	B\ = \begin{pmatrix} y & 0 \\ 0 & y^{-1}\end{pmatrix}, \\
	G = \begin{pmatrix} y^{\pm1} & 0 \\ 0 & y^{\mp1}\end{pmatrix}, 
	&\quad 
	H = \begin{pmatrix} y^{\pm1} & 0 \\ 0 & y^{\mp1}\end{pmatrix}, \quad y\in\F\setminus\{0\}.
\end{split}
\end{equation*}

The corresponding traces give curves 
\[
\begin{split}
	\cSA	&=\left\{ \left(0, y+y^{-1},  y+y^{-1},  y+y^{-1}, 0, 0, 0,  y^{1\pm1}+y^{-1\mp1}, \right. \right. \\
		&\hspace{1cm} \left. \left. y^{1\pm1}+y^{-1\mp1},  y^2+y^{-2}, 0, 0, 0,  y^{1\pm2}+y^{-1\mp2}\right)\mid  y\in\F\setminus\{0\} \right\}
\end{split} 
\]
for $l=\pm1$. 
The restrictions to $\twistKlein$ and $M$ are curves, 
\begin{align*}
	r_\klein\left(\cSA\right) 	&= \X^{\text{irr}}(\twistKlein,\F_2)=\{(0,t,0)\ \mid\ t = y+y^{-1},\ y \in\F_2\},\\
	r_M\left(\cSA\right) 		&= \X^{\text{red}}(M,\F_2)=\{(s,s,s^2)\ \mid\ s= y+y^{-1},\ y \in\F_2\},
\end{align*}
so $\cSA$ is of type \curveI.		
		
In \Cref{sec:kleinbottle}, $r_\klein\left(\cSA\right)$ is shown to detect $\annulusa$; in \Cref{sec:trefoil}, $r_M\left(\cSA\right)$ is shown to detect $\seifert$. 
We deduce both curves detect $\glueSA$ whenever $\glue=\begin{pmatrix} k & \pm1 \\ \mp1 & 0 \end{pmatrix}\in\SL_2(\Z)$.

No other curves of type \curveI\ that could give appropriate curves are found for any gluing matrix and other representations. 

Now consider the options for matrix pairs $\left(A,B\right)$ from \Cref{eqn:irredRepnKlein,eqn:redRepnKlein1,eqn:redRepnKleinChar21,eqn:redRepnKleinCharNon22} that could give appropriate curves of type \curveIII. We show there are simple closed curves in components of $N_\glue\setminus \glueSA$ that give negative valuation for each curve, proving the surface cannot be detected. 

We have already classified these options in the proof of \Cref{lem:gluingCurveTorus1,lem:gluingCurveTorus2} and only $\cTorusII$ could appear for the case \surfII.
Consider $[gh^{-1},a] \in \im(\pi_1(N_\glue\setminus \glueSA)\to\pi_1(N_\glue))$.

Using $\cTorusII$ from \Cref{eqn:cTorusII} with $\xi$ the ideal point $u\to \infty$ gives
\[
	\trFunc_{[gh^{-1},a]} = u^2 \text{ and }v_\xi(\trFunc_{[gh^{-1},a]})=-2<0, 
\]
which proves $\glueSA$ is not detected by $\cTorusII$.
This proves the result for $\glueSA$.

Consider now the options for matrix pairs $\left(A,B\right)$ from \Cref{eqn:irredRepnKlein,eqn:redRepnKlein1,eqn:redRepnKleinChar21,eqn:redRepnKleinCharNon22} in turn that could give appropriate curves of type \curveI\ that detect $\glueSB$. 

We find no ways to extend the representation from \Cref{eqn:irredRepnKlein} to get a curve of type \curveI\ that detect $\glueSB$ in any characteristic.

For general characteristic $p$ we find one way to extend the representation from \Cref{eqn:redRepnKlein1} to $N_\glue$ with a curve of type \curveI.
For gluing matrix $\glue=\begin{pmatrix} \pm1 & l \\ 0 & \mp1 \end{pmatrix}\in\SL_2(\Z)$ we get the curve of representations
\begin{equation*}
\begin{split}
	A = \begin{pmatrix} x & 1 \\ 0 & x^{-1}\end{pmatrix}, 
	&\quad 
	B = \begin{pmatrix} 1 & 0 \\ 0 & 1\end{pmatrix},\\
	G = \begin{pmatrix} x^{\pm2} & \pm(x+x^{-1}) \\ 0 & x^{\mp2}\end{pmatrix}, 
	&\quad 
	H = \begin{pmatrix} x^{\pm2} & \pm(x+x^{-1}) \\ 0 & x^{\mp2}\end{pmatrix}, \quad x\in\F\setminus\{0,1,-1\}. 
\end{split}
\end{equation*}

The corresponding traces give curves	
\begin{align*}
	\cSB =	&\left\{\left(x+x^{-1}, 2, x^2+x^{-2}, x^2+x^{-2}, x+x^{-1}, x^{1\pm2}+x^{-1\mp2}, x^{1\pm2}+x^{-1\mp2}, x^2+x^{-2}, x^2+x^{-2},  \right. \right. \\
			&\hspace{0.8cm}\left.\left. x^4+x^{-4}, x^{1\pm2}+x^{-1\mp2}, x^{1\pm2}+x^{-1\mp2}, x^{1\pm4}+x^{-1\mp4}, x^4+x^{-4}\right)\ \big|\ x\in\F\setminus\{0,1,-1\} \right\}
\end{align*}
for $k=\pm1$. 
The restrictions to $\twistKlein$ and $M$ are curves, 
\begin{align*}
	r_\klein\left(\cSB\right) 	&= \X^{\text{red}}(\twistKlein,\F_p)=\left\{\left(s,2,s\right)\ \mid\ s=x+x^{-1}, x\in\F_p\right\},\\
	r_M\left(\cSB\right) 		&= \X^{\text{red}}(M,\F_p)=\left\{\left(s,s,s^2-2\right)\ \mid\ s=x^2+x^{-2}, x\in\F_p\right\},
\end{align*}			
so $\cSB$ is of type  \curveI.

In \Cref{sec:kleinbottle}, $r_\klein\left(\cSB\right)$ is shown to detect $\annulusb$; in \Cref{sec:trefoil}, $r_M\left(\cSB\right)$ is shown to detect $\seifert$.
We deduce both curves detect $\glueSB$ whenever $\glue=\begin{pmatrix} \pm1 & l \\ 0 & \pm1\end{pmatrix}\in\SL_2(\Z)$.

No other curves of type \curveI\ that could give appropriate curves are found for any gluing matrix and other representations. 

Now consider the options for matrix pairs $\left(A,B\right)$ from \Cref{eqn:irredRepnKlein,eqn:redRepnKlein1,eqn:redRepnKleinChar21,eqn:redRepnKleinCharNon22} that could give appropriate curves of type \curveIII. 
We have already classified these options in the proof of \Cref{lem:gluingCurveTorus1,lem:gluingCurveTorus2}, which do not appear in case \surfIII.
Therefore they could not detect $\glueSB$.
This proves the result.
\end{proof}
	

\subsection{Curves in the Seifert fibered space with base orbifold $S^2(2,2,2,3)$}
\label{subsec:detection_in_seifertSphere}

The manifold $N_\glue$ forms a Seifert fibered space with base orbifold $S^2(2,2,2,3)$ if and only if $$\glue = \begin{pmatrix} k & \pm1 \\ \mp1-6k & \mp6 \end{pmatrix}\in\SL_2(\Z).$$ 
The classification of essential surfaces in $N_\glue$ in \Cref{subsec:surfaces_in_gluing} shows only the cases \surfI\ and \surfIV\ can appear. The case \surfI\ is already covered in \Cref{lem:gluingCurveTorus1,lem:gluingCurveTorus2}. The remaining case \surfIV\ and related curves in the character variety are covered in the next result. 

\begin{lemma} 
\label{lem:gluingCurveSphere} 
	Let $\F$ be an algebraically closed field of characteristic $p$. 
	Consider $N_\glue$ a Seifert fibered space with base orbifold $S^2(2,2,2,3)$ and gluing matrix $\glue$ $=$ $\begin{pmatrix} k & \pm1 \\ \mp1-6k & \mp6 \end{pmatrix}\in\SL_2(\Z).$ 
	There is a 2--dimensional component $\cRA \subset \X(N_\glue,\F)$ and for each $q\in\Z$, $\cRA$ contains a curve that detects $\glueRAp$.
\end{lemma}

Recall from \Cref{lem:gluingCurveTorus1,lem:gluingCurveTorus2}, the splitting torus $\splitTorus$ is always detected in the Seifert fibered space with base orbifold $S^2(2,2,2,3)$. Then \Cref{lem:gluingCurveTorus1,lem:gluingCurveTorus2,lem:gluingCurveSphere} prove \Cref{thm:S2family}.

\begin{proof}
We follow the approach described in \Cref{subsec:detection_in_glued}.
We previously saw the existence of the 2--dimensional component $\cRA$ in \Cref{eqn:cTorusI} as part of the proof of \Cref{lem:gluingCurveTorus1}.
We had, for the gluing matrix $\glue= \begin{pmatrix} k & \pm1 \\ \mp1-6k & \mp6 \end{pmatrix}\in\SL_2(\Z),$ the 2--dimensional components given by
\begin{equation}
\begin{split}
\label{eqn:twoDimComp}
	\cRA =	&\left\{\left(0, y+y^{-1}, (-1)^k(y+y^{-1}), (-1)^k (y+y^{-1}), 0, 0, \frac{y^2 - 1 + y^{-2}}{f (y-y^{-1})^2}+f,    \right.\right.\\
			&\hspace{1cm}\left.\left.(-1)^k(y^{1\pm1}+y^{-1\mp1}),  \pm(-1)^k\frac{y^{-1\pm2}-y^{1\mp2}}{y-y^{-1}},{1},0, \frac{y^2 - 1 + y^{-2}}{f y(y-y^{-1})^2}+fy, \right. \right. \\
			&\hspace{1cm}\left.\left. (-1)^k\left(\frac{y^2 - 1 + y^{-2}}{fy^{\pm1}(y-y^{-1})^2}+fy^{\pm1}\right), \pm\frac{y^{-1\pm1}-y^{1\mp1}}{y-y^{-1}}\right)\ \Big|\ f,y\in\F\setminus\{0\}\right\}
\end{split}
\end{equation}
for $k\in\Z$ even or odd, $l=\pm1$. 
		
We previously considered the case where $y$ is constant and obtained a curve that detects the splitting torus. 
The restrictions of the component $\cRA$ to $\twistKlein$ and $M$ are curves, 
\begin{align*}
	r_\klein\left(\cRA\right) 	&= \X^{\text{irr}}(\twistKlein,\F_p)=\left\{\left(0,s,0\right)\mid s=y+y^{-1},\ y\in\F_p\right\},\\
	r_M\left(\cRA\right) 		&= \X^{\text{irr}}(M,\F_p)=\left\{\left(s,s,1\right)\mid s=(-1)^k(y+y^{-1}),\ y\in\F_p,\ k\in\Z\right\}.
\end{align*}
In \Cref{sec:kleinbottle}, $r_\klein\left(\cRA\right)$ is shown to detect $\annulusa$; in \Cref{sec:trefoil}, $r_M\left(\cRA\right)$ is shown to detect $\annulusM$.
We deduce that each curve in $\cRA$ with the property that $I_b$ is non-constant detects $\glueRApq$ for some $(q,r)$. 
		
Now let $f=t^{u}$, $y=t^{v}$ for $t\in\F\setminus\{0\}$, $u,v \in\Z$. For each pair $(u,v)$ of co-prime integers, this gives a curve $\csub(u,v)$ inside $\cRA$ with different valuations. 
The curves are given by
\begin{equation}
\begin{split}
\label{eqn:twoDimCompSub}
	\csub(u,v) =	&\left\{\left(0, t^v+t^{-v}, (-1)^k(t^v+t^{-v}), (-1)^k (t^v+t^{-v}), 0, 0, \frac{t^{2v} - 1 + t^{-{2v}}}{t^u (t^v-t^{-v})^2}+t^u,    \right.\right.\\
				&\hspace{1cm}\left.\left.(-1)^k(t^{v(1\pm1)}+t^{-v(1\pm1)}),  \pm(-1)^k\frac{t^{v(-1\pm2)}-t^{-v(-1\pm2)}}{t^v-t^{-v}},{1},0, \frac{t^{2v} - 1 + t^{-2v}}{t^{u+v}(t^v-t^{-v})^2}+t^{u+v}, \right. \right. \\
				&\hspace{1cm}\left.\left. (-1)^k\left(\frac{t^{2v} - 1 + t^{-2v}}{t^{u\pm v}(t^v-t^{-v})^2}+t^{u\pm v}\right), \pm\frac{t^{v(-1\pm1)}-t^{-v(-1\pm1)}}{t^v-t^{-v}}\right)\ \Big|\ t\in\F\setminus\{0\}\right\}.
\end{split}
\end{equation}
Consider the recurrence relation
\begin{equation}
\label{eqn:recurRA}
	\alpha_k=\alpha_{k-1}^{-1}\alpha_{k-2}\alpha_{k-1} \text{ with initial values } \alpha_0=b^{-1}a,\ \alpha_1=a.
\end{equation}
The projection of the surface $\glueRAp$ onto the base orbifold $S^2(2,2,2,3)$ is the curve $ghg \alpha_{q}$ (see \Cref{fig:RAfamily} for examples that give an idea of the pattern). For $\glueRAp$ to be detected at the ideal point $\xi$, $ghg \alpha_{q}$ must have non-negative valuation at $\xi$. 
Consider the curve $\csub(u,v)$ with ideal point $\xi$, $t+t^{-1}\to\infty$. 
We have 
\[
\begin{split}
	\trFunc_{ghg\alpha_{q}} &= t^u t^{v(-q+1)} + \frac{t^{v(q-1)}}{t^u}\left(\frac{t^{-2v} - 1 + t^{2v}}{t^{2v}-2+t^{-2v}}\right)  \\
	\text{ and }
	v_\xi(\trFunc_{ghg\alpha_{q}}) &\geq 0 \text{ if } u = v(q-1).
\end{split}
\]
For example, $v=1$, $u=q-1$ suffices. 
Therefore the curve $\csub(q-1,1)$ in $\cRA$ defined by taking $f =t^{q-1}$, $y=t$, $t\in\F$, detects $\glueRAp$.
\end{proof}

We conjecture that we can extend this result to detect all of the surfaces $\glueRApq$ using this 2--dimensional component and have made some additional calculations in this direction.
\begin{conjecture}
\label{conj:pqSurf}
	Let $\F$ be an algebraically closed field of characteristic $p$. 
	Consider $N_\glue$ a Seifert fibered space with base orbifold $S^2(2,2,2,3)$ and gluing matrix $\glue$ $=$ $\begin{pmatrix} k & \pm1 \\ \mp1-6k & \mp6 \end{pmatrix}\in\SL_2(\Z).$ There is a curve in $\X(N_\glue,\F)$ that detects $\glueRApq$. 
Specifically, $\glueRApq$ is detected by $\csub(q-r,r)$ given in \Cref{eqn:twoDimCompSub}. 
\end{conjecture}

\begin{remark}
We have calculated several (coprime) families $(q,r)$ where we can prove \Cref{conj:pqSurf}, including: $(2a+1,2)$, $(3a+1,3)$, $(3a+2,3)$. 
This is done using similar techniques as in the proof of \Cref{lem:gluingCurveSphere}.
\end{remark}


\subsection{Curves in the Seifert fibered space with base orbifold $\mathbb{RP}^2(2,3)$}
\label{subsec:detection_in_seiferRealProj}

The manifold $N_\glue$ forms a Seifert fibered space with base orbifold $\mathbb{RP}^2(2,3)$ if and only if $$\glue = \begin{pmatrix} \pm1 & l \\ \mp6 & \pm1-6l \end{pmatrix}\in\SL_2(\Z).$$ 
The classification of essential surfaces in $N_\glue$ in \Cref{subsec:surfaces_in_gluing} shows only the cases \surfI\ and \surfV\ can appear. The case \surfI\ is already covered in \Cref{lem:gluingCurveTorus1,lem:gluingCurveTorus2}. The remaining case \surfV\ and related curves in the character variety are covered in the next result. 

\begin{lemma} 
\label{lem:gluingCurveRealProj} 
	Let $\F$ be an algebraically closed field of characteristic $p$. 
	Consider $N_\glue$ a Seifert fibered space with base orbifold $\mathbb{RP}^2(2,3)$ and gluing matrix $\glue$ $=$ $\begin{pmatrix} \pm1 & l \\ \mp6 & \pm1-6l \end{pmatrix}$ $\in$ $\SL_2(\Z)$. 
	There is a curve in $\X(N_\glue,\F)$ that detects $\glueRB$ if $p=2$ and no curve that detects $\glueRB$ if $p\neq2$.
\end{lemma}

Recall from \Cref{lem:gluingCurveTorus1,lem:gluingCurveTorus2}, the splitting torus $\splitTorus$ is never detected in the Seifert fibered space with base orbifold $\mathbb{RP}^2(2,3)$. Thus, \Cref{lem:gluingCurveTorus1,lem:gluingCurveTorus2,lem:gluingCurveRealProj} give the statement of \Cref{thm:RP2family}.

\begin{proof}
We follow the approach described in \Cref{subsec:detection_in_glued}.
The essential surface $\glueRB$ can only be detected by either a curve of type \curveI\ where the restriction to the variety of characters of $\twistKlein$ detects $\annulusM$ and the restriction to the variety of characters of $M$ detects $\annulusb$ or a curve of type \curveIII\ that has non-negative valuation on all simple closed curves in components of $N_\glue\setminus \glueRB$.

Consider first the options for matrix pairs $\left(A,B\right)$ from \Cref{eqn:irredRepnKlein,eqn:redRepnKlein1,eqn:redRepnKleinChar21,eqn:redRepnKleinCharNon22} in turn that could give appropriate curves of type \curveI. 

We find no ways to extend the representation from \Cref{eqn:irredRepnKlein} to get a curve of type \curveI\ that detect $\glueRB$ in any characteristic.

For characteristic $2$ we find one way to extend the representation from \Cref{eqn:redRepnKlein1} to $N_\glue$ with a curve of type \curveI.
For gluing matrix $\glue=\begin{pmatrix} \pm1 & l \\ \mp6 & \pm1-6l \end{pmatrix}\in\SL_2(\Z)$ we get the curve of representations
\begin{equation*}
\begin{split}
	A = \begin{pmatrix} x & 1 \\ 0 & x^{-1}\end{pmatrix}, 
	&\quad 
	B = \begin{pmatrix} 1 & 0 \\ 0 & 1\end{pmatrix}, \quad x\in\F\setminus\{0,1\} \\
	G = \begin{pmatrix} x^{\pm2} & x+x^{-1} \\ 0 & x^{\mp2}\end{pmatrix}, 
	&\quad 
	H = \begin{pmatrix} x^{\mp2} & 0 \\ \frac{1}{x+x^{-1}} & x^{\pm2}\end{pmatrix}, \quad x\in\F\setminus\{0,1\}.
\end{split}
\end{equation*}

The corresponding traces give curves			
\begin{align*}
	\cRB=& \left\{\left(x+x^{-1}, 0, x^2+x^{-2}, x^2+x^{-2}, x+x^{-1}, x^{1\pm2}+x^{-1\mp2}, x^{1\mp2}+x^{-1\pm2}+\frac{1}{x + x^{-1}}, x^2+x^{-2}, \right. \right. \\ 
		&\hspace{0.8cm}\left. \left.x^2+x^{-2},1,x^{1\pm2}+x^{-1\mp2}, x^{1\mp2}+x^{-1\pm2}+\frac{1}{x + x^{-1}}, \frac{x^{\mp2}}{x+x^{-1}}+x^{-1}, 1\right)\ \Big|\ x\in\F\setminus\{0,1\} \right\}
\end{align*}
for $k=\pm1$. 
The restrictions to $\twistKlein$ and $M$ are curves, 
\begin{align*}
	r_\klein\left(\cRB\right) 	&= \X^{\text{red}}(\twistKlein,\F_2)=\{(s,0,s)\mid s=x+x^{-1}, x \in\F_2\},\\
	r_M\left(\cRB\right) 	&= \X^{\text{irr}}(M,\F_2)=\{(s,s,1)\mid s = x^2+x^{-2}, x \in\F_2\},
\end{align*}
so $\cRB$ is of type \curveI.
		
In \Cref{sec:kleinbottle}, $r_\klein\left(\cRB\right)$ is shown to detect $\annulusb$; in \Cref{sec:trefoil}, $r_M\left(\cRB\right)$ is shown to detect $\annulusM$.
We deduce both curves detect $\glueRB$ whenever $\glue=\begin{pmatrix} \pm1 & l \\ \mp6 & \pm1-6l \end{pmatrix}\in\SL_2(\Z)$.

No other curves of type \curveI\ that could give appropriate curves are found for any gluing matrix and other representations. 

Now consider the options for matrix pairs $\left(A,B\right)$ from \Cref{eqn:irredRepnKlein,eqn:redRepnKlein1,eqn:redRepnKleinChar21,eqn:redRepnKleinCharNon22} that could give appropriate curves of type \curveIII. 
We have already classified these options in the proof of \Cref{lem:gluingCurveTorus1,lem:gluingCurveTorus2}, which do not appear in case \surfV. Therefore they could not detect $\glueRB$.
This proves the result.
\end{proof}



\bibliography{references}
\bibliographystyle{plain}


\address{Grace S. Garden\\School of Mathematics and Statistics F07, The University of Sydney, NSW 2006, Australia\\{grace.garden@sydney.edu.au\\-----}}

\address{Benjamin Martin\\Department of Mathematics, University of Aberdeen, King’s College, Aberdeen AB24 3UE, United Kingdom\\{b.martin@abdn.ac.uk\\-----}}

\address{Stephan Tillmann\\School of Mathematics and Statistics F07, The University of Sydney, NSW 2006, Australia\\{stephan.tillmann@sydney.edu.au}}


\Addresses                                       
\end{document}